\documentclass[a4paper,11pt,
]{article}
\usepackage{latexsym,amssymb,amsmath}
\usepackage{amsfonts}
\usepackage{eufrak}

\begin{document}

\newcommand{\te}{\theta}
\newcommand{\real}{\mathbb{R}}
\newcommand{\integ}{\mathbb{Z}}
\newcommand{\eps}{\varepsilon}
\newcommand{\rr}{\bf r}
\newcommand{\cc}{\bf c}
\newcommand{\zz}{\bf z}

\newcommand{\bl}{\left \{}
\newcommand{\br}{\right \}}
\newcommand{\const}{\rm const}
\newcommand{\ct}{\cos\theta}
\newcommand{\st}{\sin\theta}
\newcommand{\cp}{\cos\phi}

\newcommand{\intt}{\int\limits}
\newcommand{\summ}{\sum\limits}

\newcommand{\ham}{\vec h} 
\newcommand{\Ham}{\vec H} 
\newcommand{\hamg}{\vec g} 
\newcommand{\vp}{\frac{\dd}{\dd p}}
\newcommand{\vq}{\frac{\dd}{\dd q}}
\newcommand{\tr}{\pitchfork} 
\newcommand{\dd}{\partial}
\newcommand{\la}{\langle}
\newcommand{\ra}{\rangle}
\newcommand{\Dcurve}{\Lambda^{\circ}}
\newcommand{\ov}{\overline}
\newcommand{\Ver}{{\rm Vert}}
\newcommand{\Hor}{{\rm Hor}}
\newcommand{\cube}{\frac{3}{2}}
\newcommand{\five}[1]{(#1)^{\frac{5}{2}}}
\newcommand{\fivev}[1] {#1^{\frac{5}{2}}}
\newcommand{\Lred}{\Lambda^{\hamg}}
\newcommand{\Mu}{\mathcal{M}}
\newcommand{\DD}{\mathcal D}
\newcommand{\GG}{\mathcal G}  
\newcommand{\Iind}{{\rm Ind}_{{\Lambda_0}}\Lambda}

\newcommand{\ppotR}[3]
{

\begin{figure}\begin{center}
~\includegraphics[height=#3truecm,angle=270]{#1.eps}\\
\caption{#2}
\label{#1}
\end{center}
\end{figure}
\noindent$\!$}

\newcommand{\ppotRR}[6]
{

\begin{figure}\begin{center}
~\includegraphics[width=#5truecm,height=#6truecm,angle=270]{#1.eps}
~\includegraphics[width=#5truecm,height=#6truecm,angle=270]{#2.eps}
~\includegraphics[width=#5truecm,height=#6truecm,angle=270]{#3.eps}\\
\caption{#4}
\label{#1}
\end{center}
\end{figure}
\noindent$\!$}

\newcommand{\ead}[2]{e^{{#1}{\rm ad}{\vec{#2}}}}
\newcommand{\Chi}{\mathcal X}
\newcommand{\Ack}{ACKNOWLEGEMENTS}
\newcommand{\A}{\mathcal A}
\newtheorem{theorem}{Theorem}
\newtheorem{ex}{Example}
\newtheorem{lemma}{Lemma}
\newtheorem{property}{Property}       
\newtheorem{defi}{Definition}       
\newtheorem{remark}{{Remark}}
\newtheorem{proposition}{Proposition}
\newcommand{\BE}{\begin {equation}}
\newcommand{\EE}{\end {equation}}
\newtheorem{corollary}{Corollary}

\title {On curvatures and focal points of 
dynamical Lagrangian distributions and their reductions by 
first integrals} 
\author {Andrej A. Agrachev \thanks{ S.I.S.S.A., Via Beirut 2-4,
34013 Trieste Italy and Steklov Mathematical Institute, 
ul.~Gubkina~8, 117966 Moscow Russia; email: 
agrachev@sissa.it}\and  Natalia N. Chtcherbakova \thanks{ S.I.S.S.A., Via Beirut 2-4,
34013 Trieste Italy;  email: 
chtch@sissa.it}\and  Igor  
Zelenko \thanks{ S.I.S.S.A., Via Beirut 2-4,
34013 Trieste Italy;  email: 
zelenko@sissa.it}}  
\date{}
\maketitle

%\centerline {SISSA-ISAS, via Beirut 2-4, Trieste, Italy} 

\begin{abstract}

Pairs (Hamiltonian system, Lagrangian distribution), called 
{\sl dynamical Lagrangian distributions}, appear naturally 
in Differential Geometry, Calculus of Variations and 
Rational Mechanics. The basic differential invariants of a 
dynamical Lagrangian distribution w.r.t. the action of the 
group of symplectomorphisms of the ambient symplectic 
manifold are {\sl the curvature operator} and {\sl the 
curvature form}. These invariants can be seen as 
generalizations of the classical curvature tensor in 
Riemannian Geometry. In particular, in terms of these 
invariants one can localize the focal points along 
extremals of the corresponding variational problems. In the 
present paper we study the behavior of the curvature 
operator, the curvature form and the focal points of a 
dynamical Lagrangian distribution after its reduction by 
arbitrary first integrals in involution. The interesting 
phenomenon is that the curvature form of so-called monotone 
increasing Lagrangian dynamical distributions, which appear 
naturally in mechanical systems, does not decrease after 
reduction. It also turns out that the set of focal points 
to the given point w.r.t. the monotone increasing dynamical 
Lagrangian distribution and the corresponding set of focal 
points w.r.t. its reduction 
 by one integral are alternating sets on the corresponding integral curve of the Hamiltonian system of
 the considered  dynamical distributions. Moreover,  the first 
 focal point corresponding to the
 reduced Lagrangian distribution  comes before any focal point related to 
the original dynamical distribution.
We illustrate our results on the classical $N$-body problem. 
 
\end{abstract}
\vskip .2in


\section{Introduction}
\indent \setcounter{equation}{0}

In the present paper smooth objects are supposed to be 
$C^\infty$. The results remain valid for the class $C^k$ 
with a finite and not large $k$ but we prefer not to 
specify minimal possible $k$.

{\bf 1.1 Dynamical Lagrangian distributions.} Let $W$ be a 
symplectic manifold with symplectic form $\sigma$. 
Lagrangian distribution ${\mathcal D}$ on $W$ is a smooth 
vector sub-bundle of the tangent bundle $TW$ such that each 
fiber $\mathcal D_\lambda$ is a Lagrangian subspace of the 
linear symplectic space $T_\lambda W$, i.e., $\dim 
D_\lambda=\frac{1}{2}\dim W$ and 
$\sigma_\lambda(v_1,v_2)=0$ for all $v_1,v_2\in {\mathcal 
D}_\lambda$. For example, as a symplectic manifold one can take 
the cotangent bundle $T^*M$ 
of a manifold $M$ with standard symplectic structure and 
as a Lagrangian distribution one can take the distribution 
$\Pi(M)$ of tangent spaces to the fibers of $T^*M$, 
namely, 
\begin{equation}
\label{modex} 
\Pi(M)_\lambda=T_\lambda 
\Bigl(T^*_{\pi(\lambda)}M\Bigr), 
\end{equation}
where $\pi:\,T^*M\to M$ is the canonical projection on the 
base manifold $M$.

Let ${\mathcal H} $ be a smooth function on $W$. Denote by 
$\vec {\mathcal H}$ the Hamiltonian vector field, corresponding to 
the function  ${\mathcal H}$: 
$d{\mathcal H}(\cdot)=\sigma(\cdot,\vec {\mathcal H})$, and by 
 $e^{t \vec {\mathcal 
H}}$ the Hamiltonian flow generated by $\vec{\mathcal 
H}$.
The 
pair $(\vec{\mathcal H},{\mathcal D})$ defines the 
one-parametric family of Lagrangian distributions 
${\mathcal D}(t)=\bigl(e^{t \vec {\mathcal 
H}}\bigr)_*{\mathcal D}$. 
The pair $(\vec{\mathcal H},{\mathcal D})$ will be called 
{\it dynamical Lagrangian distribution}. 
The point $\lambda_1=e^{t_1\vec{\mathcal H}}\lambda_0$ is 
called {\sl focal} to $\lambda_0$ w.r.t. the pair 
$(\vec{\mathcal H},{\mathcal D})$ along the integral curve 
$t\mapsto e^{t\vec{\mathcal H}}\lambda_0$ of $\vec{\mathcal 
H}$, if 
\begin{equation}
\label{conj} \bigl(e^{t_1 \vec{\mathcal 
H}}\bigr)_*{\mathcal D}_{\lambda_0}\cap {\mathcal 
D}_{\lambda_1} 
\neq 0.
\end{equation}

Dynamical Lagrangian distributions appear naturally in 
Differential Geometry, Calculus of Variations and Rational 
Mechanics. The model example can be described as follows: 
\begin{ex}
\label{varex1} {\rm On a manifold $M$ for a given smooth 
function $L:TM\mapsto \mathbb{R}$, which is convex on each 
fiber, consider the following standard problem of Calculus 
of Variation with fixed endpoints $q_0$ and $q_1$ and fixed 
time $T$: 
\begin{eqnarray}
&~&A\Bigl(q(\cdot)\Bigr)=\int_{0}^T L\bigl(q(t), \dot 
q(t)\bigr)\,dt \mapsto \min \label{prob1} \\
&~&q(0)=q_0,\quad q(T)=q_1.\label{bound1} 
\end{eqnarray}
Suppose that the Legendre transform $H: T^*M\mapsto \mathbb 
{R}$ of the function $L$, 
\begin{equation}
\label{lagrtrans}
H(p,q)=\max_{X\in T_q 
M}\Bigl(p\bigl(X\bigr)- L(q,X)\Bigr), \quad q\in M, p\in 
T_q^*M,
\end{equation}
 is well defined and smooth on $T^*M$.
We will say that the 
dynamical Lagrangian distributions $\bigl(\vec H, \Pi(M)\bigr)$ 
is {\sl associated with  the 
problem} (\ref{prob1})-(\ref{bound1})\footnote {In the 
model example the Lagrangian distributions are 
integrable. For application of one-parametric 
families of non-integrable Lagrangian distributions 
see \cite{zel}.}. The curve 
$q:[0,T]\mapsto M$, satisfying (\ref{bound1}), is an 
extremal of the problem (\ref{prob1})-(\ref{bound1}) if and 
only if there exists an integral curve $\gamma:[0,T]\mapsto 
T^*M$ of $\vec H$ such that $q(t)=\pi\bigl(\gamma(t)\bigr)$ 
for all $0\leq t\le T$. 
In this case the point $\gamma(T)$ is focal to $\gamma(0)$ 
w.r.t. the pair $(\vec H, {\mathcal D})$ 
if and only if $q_1$ is conjugate to $q_0$ along the 
extremal $q(\cdot)$ in the classical variational sense for 
the problem (\ref{prob1})-(\ref{bound1}).} 
$\Box$ 
\end{ex}

The group of symplectomorphisms of $W$ acts naturally on 
Lagrangian distribution and Hamiltonian vector fields, 
therefore it acts also on dynamical Lagrangian 
distributions. Dynamical Lagrangian distributions have 
richer geometry w.r.t. this action than just Lagrangian 
distribution. For example, all integrable Lagrangian 
distributions are locally equivalent w.r.t. the action of 
the group of symplectomorphisms of $W$, while integrable 
dynamical Lagrangian distributions have functional moduli 
w.r.t. this action. 

First note that for any two vector fields $Y$, $Z$ tangent 
to the distribution $\mathcal D$ the number 
$\sigma_\lambda\bigl([\vec {\mathcal H},Y],Z\bigr)$ depends 
only on the vectors $Y(\lambda)$, $Z(\lambda)$.\footnote 
{Here $[v_1 ,v_2]$ is the Lie bracket of the vector fields 
$v_1$ and $v_2$, $[v_1,v_2]=v_1\circ v_2-v_2\circ v_1$.} 
Therefore for a given dynamical  Lagrangian 
distribution $(\vec{\mathcal H},{\mathcal 
D})$ the following bilinear form 
$Q_\lambda^{(\vec{\mathcal H},{\mathcal D})}(\cdot,\cdot)$ 
is defined on each ${\mathcal D}_\lambda$: 
\begin{equation}
\label{quad} \forall v,w\in {\mathcal D}_\lambda:\quad 
Q_\lambda^{(\vec{\mathcal H},{\mathcal 
D})}(v,w)=\sigma_\lambda\bigl([\vec {\mathcal 
H},Y],Z\bigr),\,\, Y(\lambda)=v,Z(\lambda)=w\,. 
\end{equation}
Moreover, from the fact that all ${\mathcal D}_\lambda$ are 
Lagrangian it follows that the form 
$Q_\lambda^{(\vec{\mathcal H},{\mathcal D})}$ is symmetric. 
 
A dynamical Lagrangian distribution $(\vec {\mathcal 
H},{\mathcal D})$ is called {\sl regular}, if the quadratic 
forms $v\mapsto Q^{(\vec{\mathcal H},{\mathcal 
D})}_\lambda(v,v)$ are non-degenerated for any $\lambda$. A 
dynamical Lagrangian distributions is called {\sl monotone 
(non-decreasing or non-increasing)}, if the quadratic forms 
$v\mapsto Q^{(\vec{\mathcal H},{\mathcal D})}_\lambda(v,v)$ 
are sign-definite (non-negative or non-positive definite) 
for any $\lambda$ . The regular dynamical Lagrangian 
distributions is called {\sl monotone increasing 
(decreasing)}, if the quadratic forms $v\mapsto 
Q^{(\vec{\mathcal H},{\mathcal D})}_\lambda(v,v)$ are 
 positive (negative) definite for any $\lambda$ .

\begin{remark} 
\label{secdifr} {\rm If $W=T^*M$ and $\mathcal D=\Pi(M)$ are as in 
(\ref{modex}), then the form $v\mapsto Q^{\bigl(\vec{\mathcal 
H},\Pi(M)\bigr)}_\lambda(v,v)$ coincides with the second 
differential at $\lambda$ of the restriction ${\mathcal 
H}\bigl|_{T^*_{\pi(\lambda)}M}$ of the Hamiltonian 
${\mathcal H}$ to the fiber $T^*_{\pi(\lambda)}M$. 
Therefore in this case the dynamical Lagrangian 
distribution $\bigl(\vec{\mathcal H},\Pi(M)\bigr)$ is monotone 
increasing if and only if the restrictions of $\mathcal H$ 
on each fiber of $T^*M$ are strongly convex. Consequently the dynamical 
 Lagrangian distributions $\bigl(\vec H, \Pi(M)\bigr)$ 
associated with  the 
problem  (\ref{prob1})-(\ref{bound1})
is monotone 
increasing if and only if the restrictions of the function 
$L:TM\mapsto \real$ on each fiber of $TM$ are strongly convex.}$\Box$ 
\end{remark}

 It turns out that 
under some non-restrictive assumptions on the dynamical 
Lagrangian distribution $(\vec {\mathcal H},{\mathcal D})$ 
(in particular, if this dynamical Lagrangian distribution 
is regular) one can assign to it a special linear operator 
$R^{(\vec {\mathcal H},{\mathcal D})}_\lambda$ on each 
linear spaces ${\mathcal D}_\lambda$. This operator is 
called the {\sl curvature operator} of $(\vec {\mathcal 
H},{\mathcal D})$ at $\lambda$ and it is the basic 
differential invariant of dynamical Lagrangian distribution 
$(\vec {\mathcal H},{\mathcal D})$ w.r.t. the action of the 
group of symplectomorphisms of $W$. Moreover, the following 
bilinear form 
\begin{equation}
\label{sectcurv} r^{(\vec {\mathcal H},{\mathcal 
D})}_\lambda(v, w)=Q_\lambda^{(\vec {\mathcal H},{\mathcal 
D})}\bigl(R_\lambda^{(\vec{\mathcal H},{\mathcal D})} v, w 
\bigr),\quad v,w\in {\mathcal D}_\lambda 
\end{equation}
is symmetric. The corresponding quadratic form is called 
{\it the curvature form} of the pair $(\vec {\mathcal 
H},{\mathcal D})$. Besides, the trace of the curvature 
operator 
\begin{equation}\label{Ricci}
\rho^{(\vec {\mathcal H},{\mathcal D})}_\lambda={\rm tr}R^{(\vec {\mathcal H},{\mathcal D})}_\lambda
\end{equation}
is called {\sl the generalized Ricci curvature} of $(\vec 
{\mathcal H},{\mathcal D})$ at $\lambda$. 
All these invariants where introduced in \cite{Agr2} (see 
also \cite{Agr4} and section 2 below) and the effective 
method for their calculations is given in the recent work 
\cite{Agr15}. 
Below we present the results of these calculations on 
several important examples. In all these examples $W=T^*M$ 
for some manifold $M$ and ${\mathcal D}=\Pi(M)$, a smooth 
function $L:TM\mapsto \mathbb R$ is given, the functional 
$A\bigl(q(\cdot)\bigr)$ is as in (\ref{prob1}), and 
$H:T^*M\mapsto \mathbb R$ is as in (\ref{lagrtrans}). 

\begin{ex}
({\sl Natural mechanical system}) {\rm $M=R^n$, $W=\mathbb 
R^n\times\mathbb R^n$, 
$\sigma=\sum\limits_{i=1}^ndp_i\wedge dq_i$, ${\mathcal 
D}_{(p,q)}=(\mathbb R^n,0)$, 
$L(q,X)=\frac{1}{2}\|X\|^2-U(q)$ (in this case the function 
$A\bigl(q(\cdot)\bigr)$ is the Action functional of the natural mechanical system with potential energy $U(q)$). Then 
\begin{equation}
\label{hess0}
\forall 1\leq i,j\leq n:\quad r^{(\vec H, \mathcal 
D)}_{(p,q)}(\partial_{p_i},\partial_{p_j})
=\frac{\partial^2U}{\partial q_i\partial q_j}(q).
\end{equation}
 In 
other words, in this case the curvature operator can be 
identified with the Hessian of the potential $U$.} 
\end{ex}

\begin{ex}
({\sl Riemannian manifold}) {\rm Let a Riemannian metric 
$G$ is given on a manifold $M$ by choosing an inner product 
$G_q(\cdot,\cdot)$ on each subspaces $T_qM$ for any $q\in 
M$ smoothly w.r.t. $q$. Let $L(q, X)=\frac{1}{2}G_q(X,X)$. 
The inner product 
$G_q(\cdot,\cdot)$ defines the 
canonical isomorphism between $T_q^*M$ and $T_qM$.
For any $q\in M$ and $p\in T_q^*M $ we will denote by  $p^\uparrow$
the image of $p$ under this isomorphism, namely, the vector 
$p^\uparrow\in T_qM$, 
satisfying  
\begin{equation}
\label{iso}
p(\cdot)=G_q(p^\uparrow, \cdot)
\end{equation} 
(the operation $^\uparrow$ corresponds 
to the operation of raising of indexes in the corresponding 
coordinates of co-vectors and vectors).
Since the fibers of $T^*M$ are linear spaces, one can identify 
${\mathcal D}_\lambda$ ($= T_\lambda T_{\pi(\lambda)}^*M$) with 
$T_{\pi(\lambda)}^*M$, i.e., the operation $\,^\uparrow$ is 
defined also on each ${\mathcal D}_\lambda$ with values in 
$T_{\pi(\lambda)}M$. 
It turns out (see \cite{Agr2}) that 
\begin{equation}
\label{Riemcurv}
\forall v \in {\mathcal D}_\lambda:\quad  
\bigl(R^{(\vec H, {\mathcal D})}_\lambda v\bigr)^\uparrow=
R^\nabla\bigl(
\lambda^\uparrow,
 v^\uparrow\bigr) 
\lambda^\uparrow
\end{equation}
where $R^\nabla$ is the Riemannian curvature tensor of the 
metric $G$. The right-hand side of (\ref{Riemcurv}) appears 
in the classical Jacobi equation for Jacobi vector fields 
along the Riemannian geodesics. 
Also, $\frac{1}{n-1}{\rm tr} R_\lambda^{(\vec{\mathcal 
H},{\mathcal D})}$ is exactly the Ricci curvature 
calculated at $\lambda^\uparrow$. Besides, using 
(\ref{Riemcurv}), the Riemannian curvature tensor 
$R^\nabla$ can be recovered uniquely from the curvature 
operator $R^{(\vec H, {\mathcal D})}_\lambda$. Therefore 
studying differential invariants of the appropriate 
integrable dynamical Lagrangian 
distributions, one can obtain the classical Riemannian 
tensor.} \end{ex} 

\begin{ex}
({\sl Mechanical system on a Riemannian manifold}) {\rm Let 
$G$ be the metric of the previous example and $L(q, 
X)=\frac{1}{2}G_q(X,X)-U(q)$ (in this case the function 
$A\bigl(q(\cdot)\bigr)$ is the Action functional of 
the mechanical system  on the Riemannian manifold with potential $U(q)$).
Let the operation $^\uparrow$ be as in (\ref{iso}). Then 
the curvature operator satisfies 
\begin{equation}
\label{Mechcurv} \forall v \in {\mathcal D}_\lambda:\quad 
\bigl(R^{(\vec H, {\mathcal D})}_\lambda v\bigr)^\uparrow=
R^\nabla\bigl(
\lambda^\uparrow,
 v^\uparrow\bigr) 
\lambda^\uparrow+\nabla_{v^\uparrow}({\rm grad_G\, 
U})\bigl(\pi (\lambda)\bigr), 
\end{equation}
where ${\rm grad}_GU$ is the gradient of the function $U$ 
w.r.t. the metric G, i.e., ${\rm grad}_GU=d\,U^\uparrow$, 
and $\nabla$ is the Riemannian covariant derivative.} 
$\Box$ 
\end{ex}

\begin{remark}
\label{natmon}
{\rm According to Remark \ref{secdifr}, the dynamical Lagrangian distributions from Examples 2-4 are monotone increasing.}$\Box$
\end{remark}

The generalization of different kinds of Riemannian 
curvatures, using the notion of the curvature operator of 
dynamical Lagrangian distributions, leads 
to the generalization of several classical results of 
Riemannian geometry. In \cite{Agr2} for the given monotone 
increasing dynamical Lagrangian distribution 
$(\vec{\mathcal H},\mathcal D)$ the 
estimates of intervals between two consecutive focal points 
w.r.t. the pair $(\vec {\mathcal H}, {\mathcal D})$ along 
the integral curve $\gamma(t)=e^{t\vec{\mathcal 
H}}\lambda_0$ of $\vec {\mathcal H}$ were obtained in terms 
of the curvature form of the pair $(\vec {\mathcal H}, 
{\mathcal 
D})$. 
This result is the generalization of the classical Rauch 
Comparison Theorem in 
Riemannian geometry, 
which 
gives the lower and upper bounds of the interval between 
consecutive conjugate points along the Riemannian geodesics 
in terms of upper bound for the sectional curvatures and 
lower bound for the Ricci curvature respectively. In recent 
work \cite{Agr05} it was shown that the Hamiltonian flow, 
generated by a vector field $\vec{\mathcal H}$ on the 
compact level set of ${\mathcal H}$, is hyperbolic, if 
there exists a Lagrangian distribution ${\mathcal D}$ such 
that the dynamical Lagrangian distribution $(\vec{\mathcal 
H},{\mathcal D})$ is monotone (increasing or decreasing) 
and the curvature form of so-called reduction of this 
dynamical distribution by Hamiltonian ${\mathcal H}$ on 
this level set is negative definite. This is an analog of 
the classical theorem about hyperbolicity of geodesic flows 
of negative sectional curvature on a compact Riemannian 
manifold. 

{\bf 1.2 The reduction by the first integrals.} The subject 
of the present paper is the behavior of the curvature form 
and the focal points after the reduction of the dynamical 
Lagrangian distribution $(\vec{\mathcal H},{\mathcal D})$ 
by the arbitrary $s$ first integral $g_1,\ldots, g_s$ in 
involution of the Hamiltonian ${\mathcal H}$, i.e., $s$ 
functions on $W$ such that \begin{equation} \label{pois} 
\{\mathcal{H}, g_i\}=0,\quad \{g_i, g_j\}=0,\quad \forall 
1\leq i,j\leq s
\end{equation} 
(here $\{h,g\}$ is the Poisson bracket of the functions $h$ 
and $g$, $\{h,g\}=dg(\vec h)$).
This problem appears 
naturally in the framework of mechanical systems and 
variational problems with symmetries. 
Let ${\mathcal G}=(g_1,\ldots,g_s)$ and  
\begin{equation}
\label{redef} {\mathcal D}^ {\mathcal 
{G}}_\lambda=
\Bigl(\bigcap_{i=1}^s \ker \, d_\lambda 
g_i\Bigr)\cap{\mathcal D}_\lambda+{\rm span}\bigl(\vec 
g_1(\lambda),\ldots,\vec g_s(\lambda)\bigr). 
\end{equation}
Obviously, ${\mathcal D}^{\mathcal {G}}$ is a Lagrangian 
distribution. The pair $(\vec {\mathcal H},{\mathcal 
D}^{\vec {\mathcal {G}}})$ 
is called {\sl the reduction by the tuple $\mathcal{G}$ of 
$s$ first integrals of $\mathcal{H}$ in involution} or 
shortly {\sl the $\mathcal{G}$-reduction} of the dynamical 
Lagrangian distribution $(\vec {\mathcal H},{\mathcal D})$. 
The following example justifies the word "reduction" in the 
previous definition: 

\begin{ex}
\label{ex4} {\rm Assume that we have one first integral $g$ 
of $\mathcal {H}$ such that the Hamiltonian vector field 
$\vec g$, corresponding to the first integral $g$, 
preserves the distribution ${\mathcal D}_\lambda$, namely, 
\begin{equation}
\label{presfol} (e^{t\vec g})_*{\mathcal D}={\mathcal D}. 
\end{equation}
Fixing some value $c$ of $g$, one can define (at least 
locally) the following quotient manifold: $$ W_{g,c}= 
g^{-1}(c)/{\mathcal C}\,, $$ where ${\mathcal C}$ is the 
line foliation of the integral curves of the vector field 
$\hamg$. 
The symplectic form $\sigma$ of $W$ induces the symplectic 
form 
on a manifold $W_{g, c}$, making it symplectic too. 
Besides, if we denote by $\Phi: g^{-1}(c)\mapsto W_{g, c}$ 
the canonical projection on the quotient set, the vector 
field $\Phi_*(\vec{\mathcal H} )$ is well defined 
Hamiltonian vector field on $W_{g, c}$, because by our 
assumptions the vector fields $\vec{\mathcal H}$ and 
$\hamg$ commute. Actually we have described the standard 
reduction of the Hamiltonian systems on the level set of 
the first integral, commonly used in Mechanics. In 
addition, by (\ref{presfol}), $\Phi_*({\mathcal D}^ g)$ is 
well defined Lagrangian distribution on $W_{g, c}$. So, to 
any dynamical Lagrangian distribution $( \vec{\mathcal 
H},{\mathcal D})$ on $W$ one can associate the dynamical 
Lagrangian distribution $(\Phi_*\vec{\mathcal H} , 
\Phi_*{\mathcal D}^ g)$ on the symplectic manifold $W_{g, 
c}$ of smaller dimension. 
It turns out (see subsection 2.2 
below) that the curvature form of the $g$-reduction 
$(\vec{\mathcal H} ,\DD^g)$ at $\lambda\in g^{-1}(c)$ is 
equal to the pull-back by $\Phi$ of the curvature form of 
the dynamical Lagrangian distribution $(\Phi_*\vec{\mathcal 
H},\Phi_*\DD^g)$. So, instead of $(\vec{\mathcal H} 
,\DD^g)$ one can work with $(\Phi_*\vec{\mathcal 
H},\Phi_*\DD^g)$ on the reduced symplectic space $W_{g,c}$. 
This is the essence of the reduction on the level set of 
the first integral.}$\Box$
\end{ex}

\begin{remark}
\label{redT*} {\rm Suppose now that $W=T^*M$ for some 
manifold $M$ and ${\mathcal D}=\Pi(M)$. In this case if $g$ 
is a first integral of ${\mathcal H}$, which is 
"linear w.r.t. the impulses", 
i.e., there exists a vector field $V$ on $M$ such that 
\begin{equation}
\label{linimp} 
 g(p,q) =p\bigl(V(q)\bigr),\quad  q\in M, p\in T_q^*M,
\end{equation} 
then it  satisfies
(\ref{presfol}).   
If we denote by ${\mathcal V}$ the line foliation of 
integral curves of $V$, then the reduced symplectic space 
$W_{g,c}$ can be identified with $T^*(M/{\mathcal V})$ and 
the distribution $\Pi(M)^g$ can be identified with 
$\Pi(M/{\mathcal V})$. So, after reduction we work with the 
dynamical Lagrangian distribution 
$\bigl(\Phi_*\vec{\mathcal H},\Pi(M/{\mathcal V})\bigr)$ on 
the reduced phase space $T^*(M/{\mathcal V})$ instead of 
$\bigl(\vec{\mathcal H},\Pi(M)\bigr)$. } $\Box$ 
\end{remark}

In view of the previous example the following analog of the 
notion of the focal points along the extremal w.r.t. the 
$\mathcal{G}$-reduction of the pair $(\vec {\mathcal 
H},{\mathcal D})$ is natural: The point 
$\lambda_1=e^{t_1\vec{\mathcal H}}\lambda_0$ is called {\sl 
focal} to $\lambda_0$ w.r.t. the $\mathcal {G}$-reduction 
of the pair $(\vec{\mathcal H},{\mathcal D})$ along the 
integral curve $t\mapsto e^{t\vec{\mathcal H}}\lambda_0$ of 
$\vec{\mathcal H}$, if 
\begin{equation}
\label{conjred} \Bigl(\bigr(e^{t_1 \vec{\mathcal 
H}}\bigl)_*{\mathcal D}^{\mathcal{G}}_{\lambda_0}\cap 
{\mathcal D}_{\lambda_1}^{\mathcal{G}}\Bigr)/{\rm 
span}\bigl(\vec g_1(\lambda),\ldots,\vec g_s(\lambda)\bigr) 
\neq 0. 
\end{equation} 
In the situation, described in Example \ref{ex4}, the point 
$\lambda_1=e^{t_1\vec{\mathcal H}}\lambda_0$ is focal to 
$\lambda_0$ w.r.t. the $g$-reduction of the pair 
$(\vec{\mathcal H},{\mathcal D})$ along the curve $t\mapsto 
e^{t\vec{\mathcal H}}\lambda_0$ if and only if 
$\Phi(\lambda_1)$ is focal to $\Phi(\lambda_0)$ w.r.t. the 
pair $\bigl((\Phi)_*\vec{\mathcal H},(\Phi)_*{\mathcal 
D}^g\bigr)$ along the curve $t\mapsto 
\Phi(e^{t\vec{\mathcal H}}\lambda_0)$ in $W_{g,c}$. We 
illustrate the meaning of the focal points of the reduction 
from the variational point of view on the following two 
examples. In both examples $W=T^*M$ for some manifold $M$ 
and ${\mathcal D}=\Pi(M)$, a fiber-wise convex and smooth 
function $L:TM\mapsto \mathbb R$ is given and 
$H:T^*M\mapsto \mathbb R$ is as in (\ref{lagrtrans}). 

\begin{ex} 
\label{varend} {\rm Assume that the Hamiltonian $H$ admits 
a first integral $g$, satisfying (\ref{linimp}).
It is well known that $g$, satisfying (\ref{linimp}),  
 is the first integral of $H$ if and only if
 the flow $e^{tV}$ induces the one-parametric 
 family of fiber-wise diffeomorphisms  on $TM$, 
 which preserve the function $L$, i.e., $L\circ (e^{tV})_*=L$.
 
Let ${\mathcal V}_1:{\mathbb R}\mapsto M$ be an integral 
curve of $V$ and $a(\cdot)$ be a function on $\mathcal V_1$ 
such that $$a({\mathcal V}_1(s))= s, \quad s\in {\mathbb 
R}.$$ Fix some real $c$. Then for the given point $q_0$, 
and the time $T$ consider the following variational problem 
\begin{eqnarray}
&~&\label{prob2} %A_c\Bigl(q(\cdot)\Bigr)=
\int_{0}^T L\bigl(q(t), \dot q(t)\bigr)\,dt-c \, 
a\bigl(q(T)\bigr)\to \min,\\ 
&~&\label{varbound}q(0)=q_0,\quad q(T)={\mathcal V}_1. 
\end{eqnarray}  The curve 
$q:[0,T]\mapsto M$, satisfying (\ref{varbound}), is an 
extremal of the problem (\ref{prob2})-(\ref{varbound}) if 
and only if there exists an integral curve 
$\gamma:[0,T]\mapsto g^{-1}(c)$ of $\vec H$, such that 
$q(t)=\pi\bigl(\gamma(t)\bigr)$ for all $0\leq t\le T$. 
In this case the point $\gamma(0)$ is focal to $\gamma(T)$ 
w.r.t. the $g$-reduction of the pair $(\vec H, {\mathcal 
D})$ 
if and only if the point $q_0$ is focal to the point $q(T)$ 
along the extremal $q(\cdot)$ in the classical variational 
sense for the problem (\ref{prob2})-(\ref{varbound})}. 
$\Box$ 
\end{ex}
\begin{ex}
\label{freex} {\rm Suppose that $g=\mathcal H$. For given 
real $c$ and points $q_0$, $q_1$ consider the following 
variational problem with free terminal time 
 \begin{eqnarray}
 &~&\label{prob3} %A_c\Bigl(q(\cdot)\Bigr)=
\int_{0}^T L\bigl(q(t), \dot q(t)\bigr)\,dt-c T\to \min, 
\quad T\,\,{\rm is}\,\,{\rm free,}\\ 
&~&\label{freet}q(0)=q_0,\quad q(T)=q_1. 
\end{eqnarray}
The curve $q:[0,T]\mapsto M$, satisfying (\ref{freet}), is 
an extremal of the problem (\ref{prob3})-(\ref{freet}) if 
and only if there exists an integral curve 
$\gamma:[0,T]\mapsto H^{-1}(c)$ of $\vec H$, such that 
$q(t)=\pi\bigl(\gamma(t)\bigr)$ for all $0\leq t\le T$. 
In this case the point $\gamma(0)$ is focal to $\gamma(T)$ 
w.r.t. the $H$-reduction of the pair $(\vec H, {\mathcal 
D})$ 
if and only if the point $q_0$ is focal to the point $q_1$ 
along the extremal $q(\cdot)$ in the classical variational 
sense for the problem (\ref{prob3})-(\ref{freet}). Actually 
 the considered case can be seen as a particular 
 case of the previous example.
 For this one can pass to the extended (configuration) space 
$\overline 
 M=M\times \mathbb R$ instead of $M$ and 
 take the following function $\overline L:T\overline 
 M\mapsto\mathbb R$ instead of $L$: 
 $$\overline L(\bar q, \overline X)\stackrel{def}{=}
 L(q,\frac{X}{y})y,$$ where  $\bar q\in \overline M$ such that 
 $\bar q=(q,t)$, $q\in M$, $t\in \mathbb R$ and 
 $\overline X\in T_{\bar q}\overline M$ such that $\overline X=(X,y)$, 
 $\overline X=(X, y)$, $X\in T_qM$, $y\in T_t\mathbb R\cong 
 \real$ (it is well known that $(t, H)$ is the pair of 
 conjugate variables for function $\overline L$, so as the field $V$ one takes 
 $\frac{\partial}{\partial t}$).
  }$\Box$ 
\end{ex} 

{\bf 1.3 Description of main results.} 
 For the reduced dynamical Lagrangian distribution $(\vec {\mathcal H},
 {\mathcal D}^{\mathcal G})$ 
 one can also define the curvature operator $R^{(\vec {\mathcal H},
 {\mathcal D}^{\mathcal G})}_\lambda$ and 
 the curvature form $r^{(\vec {\mathcal H},
 {\mathcal D}^{\mathcal G})}_\lambda$ on 
each linear spaces ${\mathcal D}^{\mathcal G}_\lambda$. The 
natural problem is to find the relation between $R^{(\vec 
{\mathcal H}, 
 {\mathcal D})}_\lambda$ (or $r^{(\vec {\mathcal H},
 {\mathcal D})}_\lambda$) and their  reduced analogs  $R^{(\vec {\mathcal H},
 {\mathcal D}^{\mathcal G})}_\lambda$ (or $r^{(\vec {\mathcal H},
 {\mathcal D}^{\mathcal G})}_\lambda$) on the linear space 
$\Bigl(\displaystyle{\bigcap_{i=1}^s \ker \, d_\lambda 
g_i}\Bigr)\cap{\mathcal D}_\lambda$ (which is the 
intersection of the corresponding spaces of definition 
${\mathcal D}_\lambda$ and ${\mathcal D}^{\mathcal 
G}_\lambda$). We solve this problem in section 2 for 
regular dynamical distributions. It gives an effective and 
flexible method to compute and evaluate the curvature of 
Hamiltonian systems arising in Rational Mechanics and 
geometric variational problems. The interesting phenomenon 
is that {\it the curvature form of a monotone increasing 
Lagrangian dynamical distribution does not decrease after 
reduction}. More precisely, for such distribution the 
quadratic form $$ v\mapsto r^{(\vec {\mathcal H}, 
 {\mathcal D}^\mathcal G)}_\lambda(v,v)-r^{(\vec {\mathcal H}, 
 {\mathcal D})}_\lambda(v,v), 
\quad v\in
\Bigl(\bigcap_{i=1}^s \ker \, d_\lambda 
g_i\Bigr)\cap{\mathcal D}_\lambda$$ is always {\it 
non-negative definite of rank not greater than $s$}, where 
$s$ is the number of the first integrals in the tuple 
${\mathcal G}$. 

Further, in section 3 we show that the set of focal points 
to the given point along an integral curve w.r.t. the 
monotone increasing (or decreasing) dynamical Lagrangian 
distribution and the corresponding set of focal points 
w.r.t. its reduction by one integral are alternating sets 
on the curve and the first focal point w.r.t. the reduction 
comes before any focal point w.r.t. the original dynamical 
Lagrangian distribution. In view of Examples \ref{varend} 
and \ref{freex} this result looks natural: The reduction 
enlarge the set of admissible curves in the corresponding 
variational problems (instead of the problem with fixed 
endpoints ant terminal time one obtains the problem with 
variable endpoints or free terminal time). This justifies 
the fact that the first focal point of the reduction comes 
sooner. Besides, for the mentioned examples the last fact 
and the alternation of focal points are also a consequence 
of the classical Courant Minimax Principle, applied to the 
second variation along the reference extremal in the 
corresponding variational problems. 

In addition, we demonstrate our results on the classical 
$N$ - body problem.

\section{Curvature and reduction}
\indent \setcounter{equation}{0}
 
{\bf 2.1 Curvature operator and curvature form.}
For the construction of the curvature operator of 
the dynamical Lagrangian distribution we use the theory of 
curves in the Lagrange Grassmannian, developed in \cite{Agr2} 
and \cite{Agr3}. The curve 
\begin{equation}
t\mapsto J_\lambda(t)\stackrel{def}{=}e^{-t \vec{\mathcal 
H}}_* \bigl({\mathcal D}_{e^{t \vec{\mathcal H}}\lambda}\bigr).\label{jcurve} 
\end{equation}
is called {\sl the Jacobi curve of the curve $t\mapsto 
e^{t\vec{\mathcal H}}\lambda$ attached at the point 
$\lambda$ (w.r.t. the dynamical distribution $(\vec 
{\mathcal H},{\mathcal D})$)}. It is the curve in the 
Lagrange Grassmannian $L(T_\lambda W)$ of the linear 
symplectic space $T_\lambda W$. Actually, the Jacobi curve 
is a generalization of the space of ``Jacobi fields'' along 
the extremal of variational problem of type 
(\ref{prob1})-(\ref{bound1}). Note that if 
$\bar\lambda=e^{\bar t\vec {\mathcal H}}\lambda$ then by 
(\ref{jcurve}) we have $$J_{\bar\lambda}(t)=e^{\bar t 
\vec{\mathcal H}}_*J_ \lambda(t-\bar t).$$ In other words, 
the Jacobi curves of the same integral curve of 
$\vec{\mathcal H}$ attached at two different points of this 
curve are the same, up to symplectic transformation between 
the corresponding ambient linear symplectic spaces and the 
corresponding shift of the parameterizations. Therefore, 
any differential invariants of the Jacobi curve w.r.t. the 
action of the linear Symplectic group (in other words, any 
symplectic invariant of the curve) produces the invariant 
of the corresponding dynamical Lagrangian distributions 
w.r.t. the action of the group of symplectomorphisms of the 
ambient space $W$. 

Now, following \cite{Agr2} and \cite{Agr3}, we describe the 
construction of the curvature operator of the curve in the  
Lagrange Grassmannian. Let $\Sigma$ be $2n$-dimensional 
linear space, endowed with symplectic form $\sigma$. The 
Lagrange Grassmannian $L(\Sigma)$ is real analytic 
manifold. Note that the tangent space $T_\Lambda L(\Sigma)$ 
to the Lagrangian Grassmannian at the point $\Lambda$ can 
be naturally identified with the space of quadratic forms 
${\rm Quad}(\Lambda)$ on the linear space $\Lambda\subset 
\Sigma$. Namely, take a curve $\Lambda(t)\in L(\Sigma)$ 
with $\Lambda(0)=\Lambda$. Given some vector $l\in\Lambda$, 
take a curve $l(\cdot)$ in $W$ such that $l(t)\in 
\Lambda(t)$ for all $t$ and $l(0)=l$. Define the quadratic 
form 
\begin{equation}
\label{tangdef} 
l\mapsto \sigma(\frac{d}{dt}l(0),l). 
\end{equation}
Using the fact that the spaces $\Lambda(t)$ are Lagrangian, 
it is easy to see that this form 
depends only on $\frac{d}{dt}\Lambda(0)$. So, we have the 
map from $T_\Lambda L(\Sigma)$ to the space ${\rm 
Quad}(\Lambda)$. A simple counting of dimension shows that 
this mapping is a bijection.

\begin{remark}
\label{interp} {\rm In the sequel, depending on the 
context, we will look on the elements of $T_\Lambda 
L(\Sigma)$ not only as on the quadratic forms on 
$\Lambda(t)$, but also as on the corresponding symmetric 
bilinear forms on $\Lambda(t)$ or on the corresponding 
self-adjoint operator from $\Lambda(t)$ to $\Lambda(t)^*$} 
$\Box$ 
\end{remark} 

The curve $\Lambda(\cdot)$ in $L(\Sigma)$ is called {\sl 
regular, monotone, monotone increasing (decreasing)}, if 
its velocity $\dot \Lambda(t)$ at any point $t$ is 
respectively a non-degenerated, sign-definite, positive 
(negative) definite quadratic form on the space 
$\Lambda(t)$. 
\begin{proposition}
\label{regprop} A dynamical distribution $(\vec{\mathcal 
H},{\mathcal D})$ is regular, monotone, monotone increasing 
(decreasing) if and only if all Jacobi curves w.r.t. this 
distribution are respectively regular, monotone, monotone 
increasing (decreasing) curves 
in the corresponding 
Lagrange Grassmannians. 
\end{proposition} 
 
{\bf Proof.} 
Recall that for any two vector fields 
$\vec{\mathcal H}$ and $\ell$ in $M$ one has 
\begin{equation}
\label{addt} \frac{d}{dt}\Bigl((e^{-t \vec{\mathcal H}})_*\ell \Bigr) 
=(e^{-t\vec{\mathcal H}})_*[\vec{\mathcal H},\ell]. \end{equation}
Let 
$Q_\lambda^{(\vec{\mathcal H},{\mathcal D})}$ be as in 
(\ref{quad}). Applying this fact to the Jacobi curve 
$J_\lambda(t)$ and using (\ref{quad}), (\ref{jcurve}), and 
(\ref{tangdef}) one obtains easily that 
\begin{equation}
\label{QdotJ}
 Q_\lambda^{(\vec{\mathcal H},{\mathcal 
D})}=\dot J_\lambda(0),
\end{equation}
which implies the statement of the proposition. $\Box$

Fix some $\Lambda\in L(\Sigma)$. Define the linear mapping 
$B_\Lambda:\Sigma\mapsto\Lambda^*$ in the following way: 
for given $w\in \Sigma$ one has 
 \begin{equation}
\label{gind} B_\Lambda(w)(v)=\sigma(w,v),\quad
\forall v\in \Lambda. \end{equation} Denote by 
$\Lambda^\pitchfork$ the set of all Lagrangian subspaces of 
$\Sigma$ transversal to $\Lambda$, i.e. 
$\Lambda^\pitchfork=\{\Gamma\in 
G_m(W):\Gamma\cap\Lambda=0\}$. 
Then for any subspace $\Gamma\in \Lambda^\pitchfork$ the 
restriction $B_\Lambda\bigl|_{\Gamma} \bigr 
.:\Gamma\mapsto\Lambda^*$ is an isomorphism. 

\begin{remark}
\label{ident} {\rm In other words, any $\Gamma\in 
\Lambda^\pitchfork$ can be canonically identified with the 
dual space $\Lambda^*$.} 
\end{remark}
 
Let $I_\Gamma=\bigl(B_\Lambda\bigl|_{\Gamma} \bigr 
.\bigr)^{-1}$. Note that by construction $I_\Gamma$ is 
linear mapping from $\Lambda^*$ to $\Gamma$ and 
\begin{equation}
\label{II} I_\Gamma l-I_\Delta l\in \Lambda. 
\end{equation}
The crucial observation is that the set 
$\Lambda^\pitchfork$ can be considered as an affine space 
over the linear space ${\rm Quad} (\Lambda^*)$ of all 
quadratic forms on the space $\Lambda^*$. Indeed, one can 
define the operation of subtraction on $\Lambda^\pitchfork$ 
with values in ${\rm Quad} (\Lambda^*)$ in the following 
way: \begin{equation} \label{subtr} 
(\Gamma-\Delta)(l)=\sigma(I_\Gamma l,I_\Delta l). 
\end{equation}
It is not difficult to show that $\Lambda^\pitchfork$ 
endowed with this operation of subtraction satisfies the 
axioms of affine space. For example, let us prove that 
\begin{equation}
\label{ax} (\Gamma-\Delta)+(\Delta-\Pi)=(\Gamma-\Pi) 
\end{equation}
Indeed, using skew-symmetry of $\sigma$ and relation 
(\ref{II}), one has the following series of identities for 
any $l\in \Lambda^*$ $$\sigma(I_\Gamma l,I_\Delta 
l)+\sigma(I_\Delta l,I_\Pi l)=\sigma(I_\Gamma l-I_\Pi 
l,I_\Delta l)=\sigma(I_\Gamma l-I_\Pi l,I_\Delta l-I_\Pi 
L)+$$ $$ \sigma(I_\Gamma l-I_\Pi l,I_\Pi l)=\sigma(I_\Gamma 
l-I_\Pi l,I_\Pi l)=\sigma(I_\Gamma l,I_\Pi l),$$ which 
implies (\ref{ax}).\footnote{For slightly different 
description of the affine structure on 
$\Lambda^{\pitchfork}$ see \cite{Agr2},\cite{Agr3}, and 
also \cite{Agr15}, where a similar construction is given 
for the Grassmannian $G_n(\mathbb {R}^{2n})$ of 
half-dimensional subspaces of $\mathbb {R}^{2n}$.}

Consider now some curve $\Lambda(\cdot)$ in $L(\Sigma)$. 
Fix some parameter $\tau$. Assume that 
$\Lambda(t)\in\Lambda(\tau)^\pitchfork$ for all $t$ from a 
punctured neighborhood of $\tau$. Then we obtain the curve 
$t\mapsto\Lambda(t)\in\Lambda(\tau)^\pitchfork$ in the 
affine space $\Lambda(\tau)^\pitchfork$. 
Denote by $\Lambda_\tau(t)$ the identical embedding of 
$\Lambda(t)$ in the affine space 
$\Lambda(\tau)^\pitchfork$. 
Fixing an ``origin'' $\Delta$ in $\Lambda(\tau)^\pitchfork$ 
we obtain  
 a vector function $t\mapsto \Lambda_\tau(t)-\Delta$ with values in ${\rm 
Quad}\,(\Lambda^*)$. The curve $\Lambda(\cdot)$ is called 
{\it ample} at the point $\tau$ if the function 
$\Lambda_\tau(t)-\Delta$ has the pole at $t=\tau$ 
(obviously, this definition does not depend on the choice 
of the "origin" $\Delta$ in $\Lambda(\tau)^\pitchfork$). In 
particular, if $\Lambda(\cdot)$ is a regular curve in 
$L(\Sigma)$, then one can show without difficulties that 
the function $t\mapsto \Lambda_\tau(t)-\Delta$ has a simple 
pole at $t=\tau$ for any $\Delta\in 
\Lambda(\tau)^\pitchfork$. Therefore any regular curve in 
$L(\Sigma)$ is ample at any point. 

Suppose that the curve $\Lambda(\cdot)$ is ample at some 
point $\tau$. Using only the axioms of affine space, one 
can prove easily that there exist a unique subspace 
$\Lambda^\circ(\tau)\in\Lambda^\pitchfork$ such that the 
free term in the expansion of the function $t\mapsto 
\Lambda_\tau(t)-\Lambda^\circ(\tau)$ to the Laurent series 
at $\tau$ is equal to zero. The subspace 
$\Lambda^\circ(\tau)$ is called the {\sl derivative 
subspace of the curve $\Lambda(\cdot)$ at the point 
$\tau$}. If the curve $\Lambda(\cdot)$ is ample at any 
point, one can consider the curve $\tau\mapsto 
\Lambda^\circ(\tau)$ of the derivative subspaces. This 
curve is called {\sl derivative curve} of the curve 
$\Lambda(\cdot)$. 

Now assume that the derivative curve $\Lambda^\circ(t)$ is 
smooth at a point $\tau$. In particular, the derivative 
curve of a regular curve in the Lagrange Grassmannian is 
smooth at any point (see, for example, \cite{Agr2} and the 
coordinate representation below). In general, the 
derivative curve of the ample curve is smooth at generic 
points (the points of so-called constant weight, see 
\cite{Agr3}). As was mentioned already in Remark 
\ref{interp}, one can look on $\dot\Lambda(\tau)$ and 
$\dot\Lambda^\circ(\tau)$ as on the corresponding 
self-adjoint linear mappings: 
\begin{equation}
\label{map}\dot\Lambda(\tau):  
\Lambda(\tau)\mapsto\Lambda(\tau)^*,\quad 
\dot\Lambda^\circ(\tau): 
\Lambda^\circ(\tau)\mapsto\bigl(\Lambda^\circ(\tau)\bigr)^* 
\end{equation}
Besides, by construction $\Lambda^\circ(\tau)\in 
\Lambda(\tau)^\pitchfork$. Therefore by Remark \ref{ident} 
the following spaces can be canonically identified: 
\begin{equation}
\label{comp} \Lambda(\tau)^*\cong\Lambda^\circ(\tau), \quad 
\bigl(\Lambda^\circ(\tau)\bigr)^*\cong\Lambda(\tau) 
\end{equation}
After these identifications, the composition 
$\dot\Lambda^\circ(\tau)\circ \dot\Lambda(\tau)$ is 
well-defined linear operator on $\Lambda(t)$. 

\begin{defi} 
The linear operator 
\begin{equation}
\label{curvop}
R_\Lambda(\tau)=-\dot\Lambda^\circ(\tau)\circ\dot\Lambda(\tau)\,. 
\end{equation} 
on $\Lambda(\tau)$ is called the curvature operator of the 
curve $\Lambda(\cdot)$ at a point $\tau$. The quadratic 
form $r_{\Lambda}(\tau)$ on $\Lambda(\tau)$, defined by 
\begin{equation}
\label{quadLG} 
r_{\Lambda}(\tau)(v)=\bigl(\dot\Lambda(\tau)\circ 
R_\Lambda(\tau)v\bigr)(v),\quad v\in\Lambda(\tau) 
\end{equation}
is called the curvature form of the curve $\Lambda(\cdot)$ 
at a point $\tau$.
\end{defi}

Suppose that for all Jacobi curves w.r.t. the dynamical 
Lagrangian distribution $(\vec{\mathcal H},{\mathcal D})$ 
the curvature operator is defined. {\sl The curvature 
operator $R^{(\vec {\mathcal H},{\mathcal D})}_\lambda$ 
of the dynamical Lagrangian distribution $(\vec{\mathcal 
H},{\mathcal D})$ at a point $\lambda$} is by definition 
the curvature operator of the Jacobi curve $J_\lambda(t)$ 
at $t=0$, namely, 
\begin{equation}
\label{curvopD}
 R^{(\vec {\mathcal H},{\mathcal 
D})}_\lambda=R_{J_\lambda}(0).
\end{equation}
By construction, it is the linear operator on ${\mathcal 
D}_\lambda$. {\sl The curvature form $r^{(\vec {\mathcal 
H},{\mathcal D})}_\lambda$ of the dynamical Lagrangian 
distribution $(\vec {\mathcal H},{\mathcal D})$
 at a point 
$\lambda$} is by definition the curvature form of the 
Jacobi curve $J_\lambda(t)$ at $t=0$ (see also 
(\ref{sectcurv})). 

Now for a regular curve $\Lambda(\cdot)$ in the Lagrange 
Grassmannian $L(\Sigma)$ let us give a coordinate 
representation of the derivative curve and the curvature 
operator. One can choose a basis in $\Sigma$ such that
\begin{eqnarray}
&~&\Sigma\cong\mathbb R^n\times\mathbb 
R^n=\{(x,y):x,y\in\mathbb R^n\}\label{coord},\\ 
&~&\sigma((x_1,y_1),(x_2,y_2))=\langle 
x_1,y_2\rangle-\langle x_2,y_1\rangle,\label{sympf} 
\end{eqnarray} where $\langle\cdot,\cdot\rangle$ is the standard inner product 
in $\mathbb R^n$ (such basis is called {\sl symplectic} or 
{\sl Darboux basis}). Denote by $e_i$ the $i$th vector of the 
standard basis of $\mathbb R^n$. 

Assume also that $\Lambda(\tau)\cap\{(0, y):y\in\mathbb 
R^n\}=0$. Then for any $t$ sufficiently closed to $\tau$ 
there exits the symmetric $n\times n$ matrix $S_t$ such 
that $\Lambda(t)=\{(x, S_t x): x\in \mathbb R^n\}$. The 
matrix curve $t\mapsto S_t$ is {\sl the coordinate 
representation of the curve $\Lambda(\cdot)$ (w.r.t. the 
chosen symplectic basis in $\Sigma$)}. \begin{remark} 
\label{symm} 
 Note that from (\ref{sympf}) and the 
 fact that the subspaces $\Lambda(t)$ are 
 Lagrangian it follows that the matrices $S_t$ are symmetric.
\end{remark}
 
  The curve 
$\Lambda(\cdot)$ is regular if and only if the matrices 
$\dot S_t$ are non-degenerated. The expression of the 
derivative curve and the curvature operator of the regular 
curve $\Lambda(\cdot)$ in terms of $S_t$ is given by the 
following 
\begin{proposition} 
\label{coordprop}
 The 
derivative curve $\Lambda^\circ(\tau)$ of the regular curve 
$\Lambda(\tau)$ in $L(\Sigma)$ satisfies \begin{equation} 
\label{dercurv} \Lambda^\circ(\tau)=\{ (-\frac{1}{2}\dot 
S_\tau^{-1}\ddot S_\tau \dot S_\tau^{-1} y, 
y-\frac{1}{2}S_\tau\dot S_\tau^{-1}\ddot S_\tau \dot 
S_\tau^{-1} y),y\in\mathbb R^n\}. \end{equation} In the 
basis $\{(e_i,S_\tau e_i)\}_{i=1}^n$ of $\Lambda(\tau)$ the 
curvature operator $R_\Lambda(\tau)$ is represented by the 
following matrix 
\begin{equation} \label{sch} \mathbb 
{S}(S_t)=\frac{1}{2}\dot 
S_\tau^{-1}S_\tau^{(3)}-\frac{3}{4}(\dot S_\tau^{-1} \ddot 
S_\tau)^2.\end{equation} The curvature form 
$r_\Lambda(\tau)$ has the following matrix w.r.t. the same 
basis 
\begin{equation} \label{sch1}  
-\dot S_\tau\mathbb{S}(S_t)=-\frac{1}{2}\dot 
S_\tau^{(3)}+\frac{3}{4}\ddot S_\tau\dot S_\tau^{-1} \ddot 
S_\tau .\end{equation} 
\end{proposition} 
For the proof of (\ref{dercurv}) and of the matrix 
representation (\ref{sch}) for the curvature operator see, 
for example, \cite {Agr15}. The matrix representation 
(\ref{sch1}) of the curvature form follows directly from 
(\ref{sch}) and (\ref{quadLG}). 
\begin{remark}
{\rm If $S_t$ is a scalar function (i.e., $n=1$), then 
$\mathbb{S}(S_t)$ is just the classical Schwarzian 
derivative or Schwarzian of $S_t$. 
It is well known that for scalar functions the Schwarzian 
satisfies the following remarkable identity: 
\begin{equation}
\label{mobinv} \mathbb 
S\left(\frac{a\varphi(t)+b}{c\varphi(t)+d}\right)= \mathbb 
S\bigl(\varphi(t)\bigr) \end{equation} for any constant 
$a$, $b$, $c$, and $d$, $ad-bc\neq 0$. Note that by 
choosing another symplectic basis in $\Sigma$, we obtain a 
new coordinate representation $t\mapsto \widetilde S_t$ of 
the curve $\Lambda(\cdot)$ which is a matrix M\"{o}bious 
transformation of $S_t$, 
\begin{equation}
\label{matmob} \widetilde S_t=(C+DS_t)(A+BS_t)^{-1} 
\end{equation}
for some $n\times n$ matrix $A$, $B$, $C$, and $D$. It 
turns out that the matrix Schwarzian (\ref{sch}) is 
invariant w.r.t. matrix M\"{o}bious transformations 
(\ref{matmob}) by analogy with identity (\ref{mobinv}) (the 
only difference is that instead of identity we obtain 
similarity of corresponding matrices). This is another 
explanation for invariant meaning of the expression 
(\ref{sch}), given by Proposition \ref{coordprop}.} 
\end{remark}

The coordinate representations (\ref{sch}) and (\ref{sch1}) 
are crucial in the proof of the main theorem of this 
section (see Theorem \ref{th1LG} below).

 {\bf 2.2 Curvature operator and curvature form of reduction.} 
 Now fix some $s$ vectors $l_1,\ldots l_s$ in $\Sigma$ such that
 \begin{equation}
 \label{isot} 
 \forall 1\leq i,j\leq s:\quad \sigma(l_i,l_j)=0. 
 \end{equation}
 Denote by $\ell=(l_1,\ldots l_s)$ and  ${\rm span}\, \ell={\rm 
 span}(l_1,\ldots l_s)$
 For any 
 $\Lambda\in L(\Sigma)$ let
 \begin{equation}
 \label{gamma1}
\Lambda^\ell=\Lambda\cap ({\rm span}\, \ell)^\angle+{\rm 
span}\, \ell,\quad 
\overline{\Lambda^\ell}=\Lambda^\ell\slash{\rm span}\, 
\ell, 
\end{equation}
where $({\rm span}\,\ell)^\angle\stackrel{def}{=}\{v\in 
\Sigma:\forall 1\leq i\leq s \quad \sigma(v,l_i)=0\}$ is 
the skew-orthogonal complement of the isotropic subspace 
${\rm span}\,\ell$. Actually, $\overline{\Lambda^\ell}$ is 
a Lagrangian subspace of the symplectic space $({\rm 
span}\,\ell)^\angle\slash {\rm span}\, \ell$ (with 
symplectic form induced by $\sigma$). 

Let, as before, $\Lambda(\cdot)$ be an ample curve in 
$L(\Sigma)$.
The curve $\Lambda(\cdot)^\ell$ is called {\sl the 
reduction by the $s$-tuple $\ell$}, 
satisfying (\ref{isot}), or shortly {\sl the 
$\ell$-reduction} of the curve $\Lambda(\cdot)$. Note that 
by (\ref{gamma1}) ${\rm span}\, \ell\subset 
\Lambda(t)^\ell$ for any $t$. Therefore the curve $\Lambda 
(\cdot)^\ell$ is not ample and the constructions of the 
previous subsection cannot be applied to it directly. 
Instead, suppose that the curve 
$\overline{\Lambda(\cdot)^\ell}$   
 is ample curve in 
the Lagrange Grassmannian $L(({\rm 
span}\,\ell)^\angle\slash {\rm span}\, \ell)$. Then the 
curvature operator $R_{\overline {\Lambda^\ell}}(t)$ of 
this curve is well-defined linear operator on the space 
$\overline{\Lambda(t)^\ell}$ (at least for a generic point 
$t$). Let $\phi:\Sigma\mapsto \Sigma\slash {\rm span}\, 
\ell$ be the canonical projection on the factor-space. 

\begin{defi}
The curvature operator $R_{\Lambda^\ell}(\tau)$ of the 
$\ell$-reduction $\Lambda(\cdot)^\ell$ at a point $\tau$ is 
the linear operator on $\Lambda(\tau)^\ell$, satisfying 
\begin{equation}
\label{curvred}
 R_{\Lambda^\ell}(\tau)(v)=
 \Bigl(\phi\bigl|_{\Lambda(\tau)\cap ({\rm span}\, \ell)^\angle}\Bigr)^{-1}
 \circ R_{\overline {\Lambda^\ell}}(\tau)\circ 
 \phi (v),\quad v\in\Lambda(\tau)^\ell.
 \end{equation}
 The curvature form $r_{\Lambda^\ell}(\tau)$ of the 
$\ell$-reduction $\Lambda(\cdot)^\ell$ at a point $\tau$ is 
the quadratic form on $\Lambda(\tau)^\ell$, satisfying 
\begin{equation}
\label{curvfred} 
 r_{\Lambda^\ell}(\tau)(v)=\frac{d}{d\tau}\bigl(\Lambda(\tau)^\ell\bigr)
 \bigl(R_{\Lambda^\ell}(\tau)v, v\bigr).
 \end{equation}
\end{defi}

All these constructions are directly related to the 
reduction of dynamical distributions by a tuple 
$\mathcal{G}=(g_1,\ldots, g_s)$ of $s$ involutive first 
integrals, defined in Introduction. Indeed, the Jacobi 
curves attached at some point $\lambda$ w.r.t. the 
$\mathcal G$-reduction $(\vec{\mathcal H},{\mathcal 
D}^\mathcal G)$ of a dynamical Lagrangian distribution 
$(\vec{\mathcal H},{\mathcal D})$ are exactly $\bigl(\vec 
g_1(\lambda),\ldots,\vec g_s(\lambda)\bigr) $-reductions of 
the Jacobi curves attached at $\lambda$ w.r.t. 
$(\vec{\mathcal H},{\mathcal D})$ itself. {\sl The 
curvature operator $R^{(\vec {\mathcal H}, 
 {\mathcal D}^{\mathcal G})}_\lambda$ and the curvature form $r^{(\vec {\mathcal H},
 {\mathcal D}^{\mathcal G})}_\lambda$ at $\lambda$ of the ${\mathcal G}$-reduction 
$(\vec{\mathcal H},{\mathcal D}^{\mathcal G})$} are by 
definition the curvature operator and the curvature form of 
the Jacobi curves attached at $\lambda$ w.r.t. 
$(\vec{\mathcal H},{\mathcal D}^{\mathcal G})$. 

To justify these definitions suppose that we are in 
situation of Example \ref{ex4}, i.e. ${\mathcal H}$ admits 
one first integral $g$, satisfying (\ref{presfol}). Let a 
symplectic manifold $W_{g,c}$ and a mapping 
$\Phi:g^{-1}(c)\mapsto W_{g,c}$ be as in this example. Then 
directly from the definition it follows that 
\begin{equation}
\label{pullback} \forall v\in \DD^g:\quad r_{\lambda}^ 
{(\ham,\DD^g)}(v)=r_{\Phi(\lambda)}^{(\Phi_*\ham,\Phi_*\DD^g)} 
(\Phi_*v). \end{equation} In other words, the curvature 
form of the $g$-reduction $(\vec{\mathcal H},\DD^g)$ at 
$\lambda\in g^{-1}(c)$ is equal to the pull-back by $\Phi$ 
of the curvature form of the dynamical Lagrangian 
distribution $(\Phi_*\vec{\mathcal H} ,\Phi_*\DD^g)$, 
associated to the original dynamical Lagrangian 
distribution $(\ham,\DD)$ on the reduced symplectic space 
$W_{g,c}$. 

The natural question is what is the relation between the 
curvature forms and operators of the dynamical Lagrangian 
distribution and its reduction on the common space of their 
definition. Before answering this question in the general 
situation, let us consider the following simple example: 
\begin{ex}
\label{Kepler} {\rm {\sl (Kepler's problem)}
%\label{Kepler}
Consider a natural mechanical system on $M=\real^2$ with the potential energy
$U=-r^{-1}$, where $r$ is the distance between a moving point 
in a plane and some fixed point.
This system describes the motion of the 
center of masses of two gravitationally interacting bodies in the plane of their motion (see ~\cite{Arn}).
Let $q=(r,\varphi)$ be the polar coordinates in $\real^2$.
Then the Hamiltonian function of the problem  takes the form
\begin{equation}\label{Kepler1}
h=\frac{p_r^2}{2}+\frac{p_{\varphi}^2}{2 
r^2}-\frac{1}{r}\,, \end{equation} where $p_r$ and 
$p_{\varphi}$ are the canonical impulses conjugated to $r$ 
and $\varphi$. If $\lambda=(p,q)$, where $q\in M$, $p\in 
T_q^*M$, then $p_r(\lambda)=p\bigl(\dd_r(q)\bigr)=d_q r$, 
$p_{\varphi}(\lambda)=r^2p\bigl(\dd_{\varphi}(q)\bigr)=r^2 
d_q \varphi$. Observe that $g=p_{\varphi}$ is nothing but 
the angular momentum and from (\ref{Kepler1}) we 
immediately see that it is a first integral of the system. 
Let us compare the 
curvature forms $r_{\lambda}^{(\ham, \Pi(M))}$ and 
$r_{\lambda}^{(\ham,\Pi(M)^g)}$ on the common space $\Pi(M)\cap 
{\rm ker}\,d_ \lambda g=\mathbb {R}\partial_{p_r}$ of their 
definition. 

%Consider the dynamical distribution $(\vec h, {\mathcal 
%D})$, where $\DD$ is as in (\ref{modex}). 
First, according to (\ref{hess0}) of Example 2, the 
curvature form of $(\vec h, \Pi(M))$ is equal the 
Hessian of $U$ at $q$. In particular, it implies that 
\begin{equation}
\label{kepr}
r^{\bigl(\ham, \Pi(M)\bigr)}_{\lambda}(\dd_{p_r})=\frac{\partial^2}{\partial 
r^2} U(q)= -\frac{2}{r^3}. 
\end{equation}

Further, let $c=g(\lambda)$. Note that $g$ satisfies the 
condition (\ref{linimp}) of Remark \ref{redT*} with 
$V=r^2\partial _\varphi$. Let $W_{g,c}$ and $\Phi$ be as in 
Example \ref{ex4}. Then, following Remark \ref{redT*}, 
$W_{g,c}\cong T^*\real^+$ and the dynamical Lagrangian 
distribution $(\Phi_* \vec h, \Phi_*\mathbb 
{R}\partial_{p_r})$ is equivalent (symplectomorphic) to the 
dynamical Lagrangian distribution associated with the 
natural mechanical system with configuration space 
$\real^+$ 
 and the potential energy $$ U_a=\frac{c^2}{2 
r^2}-\frac{1}{r}\, $$ ($U_a$ is the so - called {\sl 
amended} potential energy; it comes from the following 
identity: $ h\bigl|_{g^{-1}(c)}\bigr.=\frac{p_r 
^2}{2}+U_a(r)$). Hence by (\ref{pullback}) 

%The reduction by the  integral $g$ leads to the
%Hamiltonian function on the level set $g^{-1}(c)$
%
%U_a=\frac{1}{r}-\frac{c^2}{2 r^2}\,. 
%$$
%The reduced system is again a natural mechanical system, and its curvature form w.r.t $(\ham,\real\dd_{p_r})$
%is just the second derivative of $U_a$. Comparing two curvature forms on the common subspace $\ov{\DD^g}=\real\,\dd_{p_r}$
%we see
\begin{equation}
\label{kepr1}
 r_{\lambda}^{\bigl(\ham,\Pi(M)^g\bigr)}(\dd_{p_r})=\frac{d^2}{d r^2} U_a(r)=
\frac{3 
c^2}{r^4}-\frac{2}{r^3}=r^{\bigr(\ham,\Pi(M)\bigl)}_{\lambda}(\dd_{p_r})+\frac{3 
c^2}{r^4}. 
\end{equation}
%(here we used that $r_{\lambda}^
%{(\ham,\DD^g)}\dd_{p_r}=r_{\Phi(\lambda)}^{(\Phi_*\ham,\Phi_*\DD^g)}
%(\Phi_*\dd_{p_r})$, which follows from the definition of 
%the reduced curvature form) . 
Note that from (\ref{kepr1}) it follows that on the common 
space of the definition the reduced curvature form is not 
less than the curvature form itself. We will show later 
(Corollary \ref{increas1}) that this is a general 
fact.}$\Box$ 
\end{ex}

{\bf 2.3 The change of the curvature after the reduction.} 
Now we give the relation between the curvature forms of the 
regular curve $\Lambda(\cdot)$ and its $\ell$-reduction 
$\Lambda(\cdot)^\ell$, where ,as before, 
$\ell=(l_1,\ldots,l_s)$ is the tuple of $s$ vectors, 
satisfying (\ref{isot}). First, let us introduce some 
notations. Let $B_{\Lambda(t)}:\Sigma\mapsto\Lambda^*$ be 
as in (\ref{gind}). Looking at $\dot\Lambda(t)$ as at a 
linear mapping from $\Lambda(t)$ to $\Lambda(t)^*$, denote 
by $a_i(t)$, $1\leq i\leq s$ the following vectors in 
$\Lambda(t)$: 
\begin{equation}
\label{aux} a_i(t)=\bigl(\dot\Lambda(t)\bigr)^{-1}\circ 
B_{\Lambda(t)}(l_i). \end{equation} 

Using definition of $\dot\Lambda(t)$ and $B_{\Lambda(t)}$ 
one can show that $t\mapsto a_i(t)$, $a_i(t)\in\Lambda(t)$, 
is a unique vector function such that for any $t$ one has 
$$\sigma(\dot a_i(t),v)=\sigma(l_i,v),\quad \forall 
v\in\Lambda(t),$$ or, equivalently, 
\begin{equation} \label{descript} \dot a_i(t)\equiv l_i\quad 
{\rm mod}\,\, \Lambda(t). 
\end{equation}
Finally, let $A(t)$ be the $s\times s$ matrix with the 
following entries 
\begin{equation}
\label{ups} A_{km}(t)= \sigma\bigl(l_k, a_m(t)\bigr),\quad 
1\leq k,m\leq s. \end{equation} Note that by 
(\ref{descript}) and the definitions of $\dot\Lambda(t)$ 
\begin{equation}
\label{sigal} 
 A_{km}(\tau)= \sigma\bigl(\dot a_k(t),a_m(t)\bigr)
 =\dot\Lambda(t) a_k(t)\bigl(a_m(t)\bigr), 
\end{equation} 
which implies that the matrix $A(t)$ is symmetric.
\begin{theorem}
\label{th1LG} Suppose that $\Lambda(\cdot)$ is a regular 
curve in the Lagrange Grassmannian $L(\Sigma)$ and 
$\ell=(l_1,\ldots,l_s)$ is a tuple of $s$ vectors in 
$\Sigma$ such that (\ref{isot}) holds and $\det A(\tau)\neq 
0$ for some point $\tau$. Then the curvature form 
$r_\Lambda(\tau)$ of the curve $\Lambda(\cdot)$ and the 
curvature form $r_{\Lambda^\ell}(t)$ of its 
$\ell$-reduction $\Lambda(\cdot)^l$ at the point $\tau$ 
satisfy the following identity for all 
$v\in\Lambda(\tau)\cap \bigl({\rm 
span}\,\ell\bigr)^\angle$: 
\begin{equation}
\label{irr} \quad r_{\Lambda^\ell}(\tau)(v)- 
r_\Lambda(\tau)(v)=\cfrac{3}{4}\sum_{k,m=1}^s\bigr(A(\tau)^{-1}\bigl)_{km} 
\sigma\bigl(\ddot a_k(\tau),v\bigr)\sigma\bigl(\ddot 
a_m(\tau),v\bigr),\end{equation} where 
$\bigr(A(\tau)^{-1}\bigl)_{km}$ is the $km$-entry of the 
matrix $A(\tau)^{-1}$. In addition, the curvature 
operator $R_\Lambda(\tau)$ of the curve $\Lambda(\cdot)$ 
and the curvature operator $R_{\Lambda^\ell}(\tau)$ of its 
$\ell$-reduction $\Lambda(\cdot)^l$ at the point $\tau$ 
satisfy on $\Lambda(\tau)\cap\bigl({\rm span}\,\ell\bigr) 
^\angle$ the following identity: 
 \begin{equation}
\label{iRR} 
R_{\Lambda^\ell}(\tau)- 
R_\Lambda(\tau)=\cfrac{3}{4}\sum_{k,m=1}^s 
(A(\tau)^{-1})_{km}  
B_{\Lambda(\tau)} \ddot a_m(\tau)\otimes 
\Bigr(\bigl(\dot\Lambda(\tau)\bigr)^{-1}\circ 
B_{\Lambda(\tau)} \ddot a_k(\tau)\Bigl)  \,.\end{equation}
(As usual, for a given linear functional $\xi$ and a given 
vector $v$ by $\xi\otimes v$ we denote the following rank 
$1 $ linear operator $\xi\otimes v(\cdot)=\xi(\cdot) v$.)
\end{theorem} 

{\bf Proof.} 
%of Theorem \ref{th1LG} %Let us fix some $\tau$ and 
First let us prove identity (\ref{irr}). 
%for $t=\tau$. 
As before, denote by $e_i$ the $i$th vector of the standard 
basis of $\mathbb R^n$. The condition $\det A(\tau)\neq 0$ 
it equivalent to 
\begin{equation}
\label{inters} {\rm span}\bigl(a_1(\tau),\ldots, 
a_s(\tau)\bigr)\cap \bigl({\rm span}\,\ell\bigr)^\angle=0. 
\end{equation} 
 Hence one can choose a 
 basis in $\Sigma$ such that if one coordinatizes $\Sigma$ w.r.t. this basis,  
 $\Sigma\cong\mathbb 
R^n\times\mathbb R^n$,
then the symplectic form $\sigma$ satisfies (\ref{sympf}) 
and the following relations hold
\begin{eqnarray} \label{coords} &~&\Lambda(\tau)\cap 
\bigl({\rm span}\,\ell\bigr)^\angle={\rm span}\,\bigl((e_1,0),\ldots, 
(e_{n-s},0)\bigr),\\ &~&\label{dn} l_i=(0, e_{n-s+i}),\quad 
1\leq i\leq s, 
\\ &~& \label{dn1} a_i(\tau)\in {\rm span} \bigl((e_{n-s+1},0),\ldots,
(e_n,0)\bigr),\quad 1\leq i\leq s. 
\end{eqnarray} 
 
Note that by construction $$\bigl({\rm span}\,\ell\bigr)^\angle={\rm 
span}\,\bigl((e_1,0),\ldots,(e_{n-s},0), (0, e_1),
,\ldots,(0,e_n)\bigr).$$ Therefore one can make the 
following identification: 
\begin{equation}
\label{factind} ({\rm span}\,\ell)^\angle\slash {\rm 
span}\, \ell\cong {\rm 
span}\,\bigl((e_1,0),\ldots,(e_{n-s},0), (0, e_1), 
\ldots,(0,e_{n-s})\bigr). 
\end{equation}
Since by definition $a_i(t)\in\Lambda(t)$, there exists 
$b_i(t)\in \mathbb R^n$ such that $a(t)=\bigl(b_i(t), S_t 
b_i(t)\bigr)$. Note that 
\begin{equation}
\label{dota} \dot a_i(t)=\bigl(0, \dot S_t b_i(t)\bigr)+ 
\bigl(\dot b_i(t), S_t \dot b_i(t)\bigr)\equiv \bigl(0, 
\dot S_t b_i(t)\bigr)\,\,{\rm mod}\, \Lambda(t). 
\end{equation}
This, together with (\ref{descript}) and (\ref{dn}), implies 
that $l_i=\bigl(0,\dot S_t b_i(t)\bigr)$ and then 
\begin{equation}
\label{coordal} b_i(t)=\dot S_t^{-1}e_{n-s+i}. 
\end{equation} 
On the other hand, from (\ref{dn}) and (\ref{dn1}), using 
(\ref{sympf}) one can obtain that 
\begin{equation} \label{aiups} b_i(\tau)=\sum_{j=1}^s 
\sigma\bigl(a_i(\tau),l_j\bigr)e_{n-s+j}= -\sum_{j=1}^s 
A_{ij}(\tau)e_{n-s+j} 
\end{equation} 
(in the last equality we used the symmetry of the 
matrix $A(\tau)$). 
So, from (\ref{coordal}), (\ref{aiups}) and symmetry of 
$\dot S_\tau$ (see Remark \ref{symm}) it follows that 

\begin{eqnarray} 
&~& \label{block1}\forall 1\leq i\leq n-s,\quad n-s+1\leq j\le 
n:\quad (\dot S_\tau^{-1})_{ij}=(\dot 
S_\tau^{-1})_{ji}=0,\\ &~&\label{block2} \forall 1\leq 
i,j\leq s:\quad (\dot 
S_\tau^{-1})_{n-s+i,n-s+j}=-A_{ij}(\tau). 
\end{eqnarray}

Further, for given $n\times n$ matrix ${\mathcal A}$ denote 
by $C({\mathcal A})$ the $(n-s)\times(n-s)$ matrix, 
obtained from ${\mathcal A}$ by erasing the last $s$ 
columns and rows. Consider the curve 
$\overline{\Lambda(\cdot)^\ell}$ in the Lagrange 
Grassmannian $L(({\rm span}\,\ell)^\angle\slash {\rm 
span}\, \ell)$ (see (\ref{gamma1}) for the notation). By 
construction, if $S_t$ is the coordinate representation of 
the curve $\Lambda(t)$ w.r.t. the chosen symplectic basis, 
then $C(S_t)$ is a coordinate representation of the curve 
$\overline{\Lambda(\cdot)^l}$ w.r.t. the basis of $({\rm 
span}\,\ell)^\angle\slash {\rm span}\, \ell$, indicated in 
(\ref{factind}). Hence from (\ref{block1}), (\ref{block2}), 
and assumption $\det\, A(\tau)\neq 0$ it follows that the 
germ at $\tau$ of the curve 
$\overline{\Lambda(\cdot)^\ell}$ is regular. In particular, 
the curvature form of the $\ell$-reduction 
$\Lambda(\cdot)^\ell$ is well defined at $\tau$. Using 
(\ref{sch1}), we obtain that the quadratic form 
$$r_{\Lambda^\ell}(\tau)\Bigl|_ {\Lambda(\tau)\cap ({\rm 
span}\,\ell)^\angle }\Bigr.-r_{\Lambda}(\tau)\Bigl|_ 
{\Lambda(\tau)\cap ({\rm span}\,\ell)^\angle}\Bigr.,$$ 
%defined on $\Lambda(\tau)\cap ({\rm span}\,\ell)^\angle$, 
has the following 
matrix 
in the basis $\bigl((e_1,0),\ldots,(e_{n-s},0)\bigr)$: 
\begin{equation*}
- \frac{d}{d\tau} 
C(S_\tau)\mathbb{S}\bigl(C(S_\tau)\bigr)+C\bigl( \dot 
S_\tau\mathbb{S}\bigl(S_\tau)\bigr). 
\end{equation*}
Using the blocked structure of the matrix $\dot S_\tau$, 
given by (\ref{block1}), 
one can obtain from (\ref{sch1}) without difficulties that 
\begin{equation*}
-\frac{d}{d\tau} 
C(S_\tau)\mathbb{S}\bigl(C(S_\tau)\bigr)+C\bigl( \dot 
S_\tau\mathbb{S}\bigl(S_\tau)\bigr)=-\frac{3}{4}\Bigl\{ 
\sum_{k,m=n-s+1}^n (\ddot S_\tau)_{ik}(\dot 
S_\tau^{-1})_{km} (\ddot S_\tau)_{mj}\Bigr\}_{i,j=1}^{n-s} 
\end{equation*} 
In order to prove (\ref{irr}), it is sufficient to prove 
the following 
\begin{lemma}
The restriction of the quadratic form 
\begin{equation}
\label{ddquad} v\mapsto\sum_{k,m=1}^s(A(\tau)^{-1})_{km} 
\sigma\bigl(\ddot a_k(\tau),v\bigr)\sigma\bigl(\ddot 
a_m(\tau),v\bigr) 
\end{equation}
on $\Lambda(\tau)\cap l^\angle$ has the matrix with $i 
j$-entry equal to 
\begin{equation*}
-\sum_{k,m=n-s+1}^n (\ddot S_\tau)_{ik}(\dot 
S_\tau^{-1})_{km} (\ddot S_\tau)_{mj}
\end{equation*}
in the basis $\bigl((e_1,0),\ldots,(e_{n-s},0)\bigr)$. 
\end{lemma} 

{\bf Proof.} First, using the symmetry of $S_t$ and 
(\ref{block2}), one has the following identity: 
\begin{equation}
\label{preob} \begin{split}\sum_{k,m=n-s+1}^n (\ddot 
S_\tau)_{ik}(\dot S_\tau^{-1})_{km} (\ddot S_\tau)_{mj}= 
\sum_{k,m=n-s+1}^n (\ddot S_\tau\dot S_\tau^{-1})_{ik}(\dot 
S_\tau)_{km} (\ddot S_\tau \dot S_\tau^{-1})_{jm}=\\ 
-\sum_{k,m=1}^s (\ddot S_\tau\dot 
S_\tau^{-1})_{i,n-s+k}\bigl(A(\tau)^{-1}\bigr)_{km} (\ddot 
S_\tau \dot S_\tau^{-1})_{j,n-s+m}\end{split} 
\end{equation} 
On the other hand, from (\ref{dota}) we have: $$\ddot 
a_i(t)=\bigl(0, \ddot S_t b_i(t)\bigr)+ 2\bigl(0, \dot S_t 
\dot b_i(t)\bigr) +\bigl(\ddot b_i(t), S_t \ddot 
b_i(t)\bigr).$$ 

Substitute $t=\tau$ in the last relation. Note that $\dot 
S_\tau$ has the same blocked structure, as $\dot 
S_\tau^{-1}$, which together with (\ref{aiups}) implies 
that $\bigl(0, \dot S_t \dot b_i(t)\bigr)\in {\rm 
span}\,\ell$. From this and (\ref{coordal}) it follows that 
\begin{equation}
\label{ddota} \ddot a_i(\tau) \equiv \bigl(0, \ddot S_\tau 
\dot S_\tau^{-1} e_{n-s+i}(t)\bigr)\quad {\rm mod}\,\, {\rm 
span} \bigl(\Lambda(\tau),l_1,\ldots,l_s\bigr). 
\end{equation}
 Hence 
\begin{equation}
\label{fin} \forall 1\leq j \leq n-s, 1\leq i\leq s:\quad 
\sigma \bigl (\ddot a_i(\tau),(e_j,0)\bigr)= 
-(\ddot S_\tau \dot S_\tau^{-1})_{j,n-s+i}, 
\end{equation} 
The last identity together with (\ref{preob}) implies the 
statement of the lemma 
and also formula (\ref{irr}). $\Box$ 

Finally, identity (\ref{iRR}) follows directly from 
(\ref{quadLG}). The proof of the theorem is completed. 
$\Box$ \vskip .2in

Note that if the curve $\Lambda(\cdot)$ is monotone 
increasing or decreasing , then from (\ref{sigal}) it 
follows that the condition $\det A(\tau)\neq 0$ is 
equivalent to the following condition 
\begin{equation}
\label{inters1} \Lambda(\tau)\cap {\rm span}\, \ell=0. 
\end{equation} 
As a direct consequence of identity (\ref{irr}) and the 
last fact one has the following 
\begin{corollary} \label{increas} If the curve 
$\Lambda(\cdot)$ is monotone increasing and the tuple 
$\ell=(l_1,\ldots,l_s)$ of $s$ vectors in $\Sigma$ 
satisfies (\ref{inters1}) for some $t$, then on the space 
$\Lambda(t)\cap \bigl({\rm span}\, \ell\bigr)^{\angle}$ the 
curvature form of the $\ell$-reduction of the curve 
$\Lambda(\cdot)$ is not less than the curvature form of the 
curve $\Lambda(\cdot)$ itself. Moreover, on the space 
$\Lambda(t)\cap \bigl({\rm span}\, \ell\bigr)^{\angle}$ the 
difference between the curvature form of the 
$\ell$-reduction of the curve $\Lambda(\cdot)$ and the 
curvature form of the curve $\Lambda(\cdot)$ itself is 
non-negative definite quadratic form of rank not greater 
than $s$. 
\end{corollary} 

Let us translate the results of Theorem \ref{th1LG} in 
terms of a regular dynamical Lagrangian distribution 
$(\vec{\mathcal H},{\mathcal D})$ on a symplectic space 
$W$. Suppose that the Hamiltonian ${\mathcal H}$ admits a 
$s$-tuple $\mathcal{G}=(g_1,\ldots, g_s)$ of involutive 
first integrals. 
Similarly to 
(\ref{gind}) denote by $\mathcal{B}_{\mathcal 
{D}_\lambda}:T_\lambda W\mapsto \mathcal{D}_\lambda^*$ the 
linear mapping such that for given $Y\in T_\lambda W$ the 
following identity holds 
\begin{equation}
\label{gind1} \mathcal{B}_{\mathcal 
{D}_\lambda}Y(Z)=\sigma(Y,Z),\quad 
\forall Z\in \mathcal {D}_\lambda. \end{equation} Let us 
look at $\dot J_\lambda(0)$ (the velocity at $t=0$ of the 
Jacobi curve attached at $\lambda$) as at a linear mapping 
from ${\mathcal D}_\lambda$ to ${\mathcal D}_\lambda^*$. 
Then using regularity  by analogy with 
(\ref{aux}) one can define the following  $s$ vector fields ${\mathcal X}_i$ on 
$W$: 
\begin{equation}
\label{aux1} \mathcal{X}_i(\lambda)=\bigl(\dot 
J_\lambda(0)\bigr)^{-1}\circ \mathcal{B}_{\mathcal 
{D}_\lambda}\bigl(\vec g_i(\lambda)\bigr). 
\end{equation} 
Using relation (\ref{addt}) one can obtain by analogy with 
(\ref{descript}) that $\mathcal{X}_i$ is a unique vector 
field, satisfying $\mathcal {X}_i(\lambda)\in 
\mathcal{D}_\lambda$ and 
\begin{equation}
\label{derX} \quad [\vec{\mathcal H},{\mathcal 
X}_i](\lambda)\equiv \vec g_i(\lambda)\quad {\rm 
mod}\,\,{\mathcal D}_\lambda
\end{equation}
for all $\lambda\in W$. Finally let $\Upsilon(\lambda)$ be 
the $s\times s$ matrix with the following entries 
\begin{equation}
\label{ups1} \Upsilon(\lambda)_{km}= 
\sigma_\lambda\bigl(\vec g_k, {\mathcal X}_m\bigr),\quad 
1\leq k,m\leq s. 
\end{equation} 

From Theorem \ref{th1LG} and relation (\ref{addt}) one gets 
immediately 
\begin{theorem}
\label{th2LG} Suppose that $(\vec{\mathcal H},{\mathcal 
D})$ is a regular Lagrangian dynamical distribution on a 
symplectic space $W$ and the Hamiltonian ${\mathcal H}$ 
admits a tuple $\mathcal{G}=(g_1,\ldots, g_s)$ of $s$ 
involutive first integrals such that 
$\det\,\Upsilon(\lambda)\neq 0$. Then the curvature form 
$r^{(\vec {\mathcal H},{\mathcal D})}_\lambda$ of the 
dynamical distribution $(\vec{\mathcal H},{\mathcal D})$ 
and the curvature form $r^{(\vec {\mathcal H},{\mathcal 
D^{\mathcal G}})}_\lambda$ of its ${\mathcal G}$-reduction 
$(\vec {\mathcal H},{\mathcal D^{\mathcal G}})$ satisfy the 
following identity for all $v\in 
\Bigl(\displaystyle{\bigcap_{i=1}^s \ker \, d_\lambda 
g_i}\Bigr)\cap{\mathcal D}_\lambda$ 
\begin{equation}
\label{irr1} \quad r^{(\vec {\mathcal H},{\mathcal 
D^{\mathcal G}})}_\lambda(v)- r^{(\vec {\mathcal H},{\mathcal 
D})}_\lambda(v)=\cfrac{3}{4}\sum_{k,m=1}^s 
\bigl(\Upsilon(\lambda)^{-1}\bigr)_{km} 
\sigma_\lambda\bigl(\bigl[\vec{\mathcal H}, 
[\vec{\mathcal H},\mathcal{X}_k]\bigr], v 
\bigr)\sigma_\lambda\bigl(\bigl[\vec{\mathcal H}, 
[\vec{\mathcal H},\mathcal{X}_m]\bigr], v\bigr) 
,\end{equation} while the curvature operator $R^{(\vec 
{\mathcal H},{\mathcal D})}_\lambda$ of the dynamical 
distribution $(\vec{\mathcal H},{\mathcal D})$ and the 
curvature operator $R^{(\vec {\mathcal H},{\mathcal 
D}^{\mathcal G})}_\lambda$ of its $\mathcal G$-reduction 
$(\vec {\mathcal H},{\mathcal D}^{\mathcal G})$ satisfy on 
$\Bigl(\displaystyle{\bigcap_{i=1}^s \ker \, d_\lambda 
g_i}\Bigr)\cap{\mathcal D}_\lambda$ the following identity: 
 \begin{equation}
\label{iRR1} 
\begin{split}
&\hspace {1.5in} R^{(\vec {\mathcal H},{\mathcal 
D}^{\mathcal G})}_\lambda-R^{(\vec {\mathcal H},{\mathcal 
D})}_\lambda =\\ 
&\cfrac{3}{4}\sum_{k,m=1}^s\bigl(\Upsilon(\lambda)^{-1}\bigr)_{km} 
\mathcal{B}_{\mathcal {D}_\lambda} \bigl[\vec{\mathcal H}, 
[\vec{\mathcal H},\mathcal{X}_m]\bigr](\lambda)\otimes 
\Bigr(\bigl(\dot J_\lambda(0)\bigr)^{-1}\circ 
\mathcal{B}_{\mathcal {D}_\lambda} \bigl[\vec{\mathcal H}, 
[\vec{\mathcal H},\mathcal{X}_k]\bigr](\lambda) \Bigl) 
.
\end{split}
\end{equation} 
\end{theorem}

Also, by analogy with Corollary \ref{increas} we have

\begin{corollary} \label{increas1} If the dynamical Lagrangian 
distribution $(\vec{\mathcal H},{\mathcal D})$ is monotone 
increasing and the Hamiltonian ${\mathcal H}$ admits a 
tuple $\mathcal{G}=(g_1,\ldots, g_s)$ of $s$ involutive 
first integrals such that
\begin{equation}
\label{inters2} D_\lambda\cap{\rm span} \bigl(\vec 
g_1(\lambda),\ldots, \vec g_s(\lambda)\bigr)=0,
\end{equation}
then on the space $
\Bigl(\displaystyle{\bigcap_{i=1}^s \ker \, d_\lambda 
g_i}\Bigr)\cap{\mathcal D}_\lambda$ the curvature form of 
the $\mathcal{G}$-reduction of the dynamical Lagrangian 
distribution $(\vec{\mathcal H},{\mathcal D})$ is not less 
than the curvature form of $(\vec{\mathcal H},{\mathcal 
D})$ itself. 
 Moreover, on the space $%\ov{\DD^{\GG}_{\lambda}}
 \Bigl(\displaystyle{\bigcap_{i=1}^s \ker \, d_\lambda 
g_i}\Bigr)\cap{\mathcal D}_\lambda 
 $ the difference 
$$r^{(\vec {\mathcal H},{\mathcal 
D^{\mathcal G}})}_\lambda- r^{(\vec {\mathcal H},{\mathcal 
D})}_\lambda$$
is non-negative definite quadratic form of rank not 
greater than $s$. 
\end{corollary}

Now let us give the coordinate representation of the vector 
fields ${\mathcal X}_i$, $1\leq i\leq s$ from Theorem 
\ref{th2LG} in the case, when $W=T^*M$ and 
${\mathcal D}=\Pi(M)$. Let $q=(q^1,\ldots,q^n)$ be local 
coordinates in some open subset ${\mathcal N} $ of $M$ and 
$p=(p_1,\ldots,p_n)$ be induced coordinates in the fiber of 
$T^*{\mathcal N}$ so that the canonical symplectic form is given by 
$\sigma=\sum\limits_{i=1}^ndp_i\wedge dq^i$. It gives the 
identification of $T^*{\mathcal N}\cong \mathbb R^n\times \mathbb 
R^n=\{(p,q),p,q\in \mathbb R^n)$ (so, ${\mathcal N}=0\times \mathbb R^n$). 
Also the tangent space $T_\lambda (T^*{\mathcal N})$ to $T^*{\mathcal N}$ 
at any 
$\lambda$ is identified with $ \mathbb R^n\times \mathbb 
R^n$. Under this identification $\vec{\mathcal 
H}=\left(-\frac{\partial{\mathcal H}}{\partial 
q},\frac{\partial{\mathcal H}}{\partial p}\right)$, where for 
given function $h$ on $T^*{\mathcal N}$ we denote by $\frac{\partial 
h}{\partial q}=\left(\frac{\partial h}{\partial 
q^1},\ldots,\frac{\partial h}{\partial q^n}\right)^T$ and 
$\frac{\partial h}{\partial p}=\left(\frac{\partial h}{\partial 
p_1},\ldots,\frac{\partial h}{\partial p_n}\right)^T$. Denote  
by ${\mathcal H}_{pp}$ the Hessian matrix of the 
restriction of ${\mathcal H}$ to the fibers. Then from 
Remark \ref{secdifr} and relation (\ref{aux1}) we have 
\begin{equation}
\label{en} {\mathcal X}_i=({\mathcal H}_{pp}^{-1}\frac{\partial g_i}{\partial 
p},0) 
\end{equation} 

Now suppose for simplicity that the dynamical Lagrangian 
distribution is associated with a natural mechanical system 
(Example 2) or, more generally, with a mechanical system on 
a Riemannian manifolds (Example 3), which admits one or 
several first integrals being in involution and linear 
w.r.t. the impulses. One way to compute the reduced 
curvatures is to pass to the reduced phase space, as was 
described in Remark \ref{redT*}, and apply the method of 
computation of the curvatures from \cite{Agr15} to the 
corresponding dynamical Lagrangian distribution in the 
reduced phase space (this way was actually implemented in 
Example \ref{Kepler}). But in order to apply the method of 
\cite{Agr15} we need to find a new canonical coordinates in 
the reduced phase space, which is not just a trivial 
exercise. Moreover, very often the new Hamiltonian system 
on the reduced phase space has more complicated form than 
the original one. Both these facts make the computation in 
this way quite tricky. Theorem 2 gives another method to 
compute all reduced curvatures without passing to the 
reduced phase space: to do this one can combine 
(\ref{hess0}) or (\ref{Mechcurv}) with (\ref{irr1}) (or 
(\ref{iRR1})) and (\ref {en}). This method is more 
effective from the computational point of view, especially 
if the number of the involutive first integrals is 
essentially less than the number of the degrees of freedom 
in the problem. We illustrate the effectiveness of this 
method on the following example: 

\begin{ex}{\sl( Plane $N$-body problem with equal masses)} {\rm
Let us consider the motion 
of $N$ bodies of unit mass in $\real^2$ endowed with the
standard Cartesian coordinates  so that  
$r_i=(q_{2i-1},q_{2i})\in \real^2$ represents the radius 
vector of the $i$-th body with respect to some inertial 
frame. It is described by a natural mechanical system on $M=\real^{2N}$ with potential energy
\begin{equation}
U(r_1,\dots,r_N)=-\summ_{i<j}^N 
\frac{1}{r_{ij}}\,,\qquad r_{ij}=\|r_i-r_j\|\,.\label{potential} \end{equation}

Then 
$T^*M\cong\real^{2N}\times\real^{2N}=\{(p,q),\;p,q\in\real^{2N}\}$, 
$p_1\dots,p_{2N}$ are the canonical impulses conjugated to 
$q_1,\dots,q_{2N}$ ($p_i\sim \dot q_i$). The systems has 
the following first integral 
\begin{equation}
\label{intg}
 g=\summ_{i=1}^N 
(p_{2i}q_{2i-1}-p_{2i-1}q_{2i})\, 
\end{equation}
 which is nothing but 
the angular momentum (in the considered planar case the 
angular momentum is scalar). From Example 2 we know that 
the generalized curvature form of 
 the dynamical Lagrangian distribution $(\Ham,\Pi(M))$
is just the Hessian of the potential energy $U$
and the generalized Ricci curvature (see (\ref{Ricci}) for the definition) is the Laplacian of $U$, which can be calculated without difficulties:
\begin{equation}
\label{laplacian}
\rho_{\lambda}^{\bigl(\Ham,\Pi(M)\bigr)}=\Delta\, U=-2\summ_{i<j}^N \frac{1}{r_{ij}^3}, \quad \lambda=(p,q).
\end{equation}

Our goal is to compute the reduced generalized Ricci 
curvatures $\rho_{\lambda}^{\bigl(\Ham,\Pi(M)^g\bigr)}$ , 
using the formula (\ref{irr1}). In our case $s=1$. Let 
${\mathcal X}$ be as in (\ref{aux1}) with $g$ instead of 
$g_i$. Note that by definition the vector ${\mathcal 
X}(\lambda)$ is orthogonal to the subspace 
$\Pi(M)_{\lambda} \cap {\rm ker}d_{\lambda}g$ w.r.t. the 
inner product 
$Q_\lambda^{\bigr(\Ham,\Pi(M)\bigl)}(\cdot,\cdot)$. 
%As was mentioned before, by analogy with the classical Rauch theorem 
%of Riemannian the estimation from below of the generalized Ricci curvature 
%gives the estimation from above of successive focal points along the integral 
%curve of $\vec H$ w.r.t. the 
%corresponding dynamical Lagrangian distribution.  
Therefore 
\begin{equation}
\label{trace}
\begin{split} 
&\rho_{\lambda}^{\bigl(\Ham,\Pi(M)^g\bigr)}=\rho_{\lambda}^{\bigr(\Ham,\Pi(M)\bigl)}-
\cfrac{r_{\lambda}^{\bigr(\Ham,\Pi(M)\bigl)}({\mathcal X})}
{Q_{\lambda}^{\bigr(\Ham,\Pi(M)\bigl)}({\mathcal X},{\mathcal X})} +\\
&{\rm tr}\left[\left(R^{\bigr(\Ham,\Pi(M)^{g}\bigl)}_\lambda- 
R^{\bigl(\Ham,\Pi(M)\bigr)}_\lambda\right)\Big|_{\Pi(M)_{\lambda}
\cap {\rm ker}d_{\lambda}g}\right]
%\,\quad\lambda=(p,q).
\end{split}
\end{equation}
Further, the Hamiltonian vector fields corresponding to the functions  $H$ and $g$ are given by 
$\Ham=(U_q,p)$ and $\hamg=(Jp,-Jq)$ with $U_q=\left(\frac{\partial U}{\partial 
q_1},\ldots,\frac{\partial U}{\partial q_n}\right)^T$ and
$J$ being the unit symplectic $2N\times 2N$ matrix:
\begin{equation*}
J=\left(\begin{array}{ccccc}
0&1&0&\dots&0\\
-1&0&0&\dots&0\\
\;&\;&\ddots&\;&\;\\
0&\dots&0&0&1\\
0&\dots&0&-1&0
 \end{array}\right)\,.
%\label{iSdot0}
\end{equation*}
Applying  formula (\ref{en}) we find that
\begin{equation}
\label{chiJ}
{\mathcal X}=(Jq,0).
\end{equation}
Denote $\ov{\mathcal X}=Jq$. Also, let $\la\cdot, \cdot\ra$ and $\|\cdot\|$ be 
the standard Euclidean inner product and norm. 
Using again (\ref{hess0}), one can obtain by direct computation that
%\rho_{\lambda}^{(\Ham,\DD)}=\Delta U(q)\,,\qquad
\begin{equation}
\label{secterm}
r_{\lambda}^{(\Ham,\Pi(M))}({\mathcal X})=\la U_{qq} \ov{\mathcal X},\ov{\mathcal X}\ra=-U.
\end{equation}
Further, using (\ref{irr1}) (or (\ref{iRR1})) and (\ref{chiJ}), 
one can obtain without 
difficulties that  
\begin{equation}
\label{irrf}
\begin{split}
&{\rm 
tr}\left[\left(R^{\bigl(\Ham,\Pi(M)^{g}\bigr)}_\lambda- 
R^{\bigl(\Ham,\Pi(M)\bigr)}_\lambda\right)\Big|_{\Pi(M)_ 
{\lambda}\cap {\rm ker}d_{\lambda}g}\right]=\\ 
&\frac{3}{4\sigma({\mathcal 
X},\hamg)}\left(\summ_{i=1}^{2N} 
\sigma([\Ham,[\Ham,{\mathcal X}]],\dd _{p_i})^2- 
\frac{\sigma([\Ham,[\Ham,{\mathcal X}]],{\mathcal X})^2} 
{Q_{\lambda}^{\bigr(\Ham,\Pi(M)\bigl)}({\mathcal 
X},{\mathcal X})}\right)=\frac{3}{4}\left(\frac{\| 
p\|^2}{\|q\|^2}- \frac{\la p,q\ra^2}{\|q\|^4}\right). 
\end{split}
\end{equation}
Substituting (\ref{laplacian}), (\ref{secterm}), and (\ref{irrf}) into 
(\ref{trace}) we obtain finally that
\begin{equation}
\label{Nricci}
\rho_{\lambda}^{\bigl(\Ham,\Pi(M)^g\bigr)}
%=\Delta U-\frac{\la U_{qq} \ov{\mathcal X},\ov{\mathcal X}\ra}{\|\ov{\mathcal X}\|^2}+=
%$$
%\BE\label{Nricci}
=-2\summ_{i<j}^N \frac{1}{r_{ij}^3}-\frac{U}{I}+
\frac{3}{I^2}(2 T I-\frac{1}{4}\{H, I\}^2)\,,
\end{equation}
where $I=\|q\|^2$, $T=\frac{1}{2}\|p\|^2$ are the central 
momentum of inertia and the kinetic energy of the system of 
$N$ bodies. Note that the sum of the first two terms in 
(\ref{Nricci}) is the trace of the restriction of the 
curvature operator $R^{\bigl(\Ham,\Pi(M)\bigr)}_\lambda$ on 
the space $\Pi(M)_ {\lambda}\cap {\rm ker}d_{\lambda}g$. 
So, by Remark \ref {natmon} and Corollary \ref{increas1}, 
the last term in (\ref{Nricci}) has to be nonnegative. 
Actually this term contains the right-hand side of the 
famous Sundman's inequality $2 T I-\frac{1}{4}\{H, I\}^2\ge 
0$ and it is nothing but the generalized area of the 
parallelogram formed by two $2N$-dimensional vectors $p$ 
and $q$. } $\Box$ 
\end{ex}

\section{Focal points and Reduction}
\indent
\setcounter{equation}{0}

In the present section we study the relation between the 
set of focal points to the given point w.r.t. the monotone 
increasing (or decreasing) dynamical Lagrangian 
distribution and the set of focal points 
w.r.t. its reduction. As before, first we prove the 
corresponding result for the curves in Lagrange 
Grassmannians and then reformulate it in terns of the 
dynamical Lagrangian distributions. 

Let $\Lambda(\cdot)$ be a curve in the Lagrange Grassmannian 
$L(\Sigma)$, defined on the interval $[0, T]$. The time 
$t_1$ is called {\sl focal to the time $0$ w.r.t. the curve 
$\Lambda(\cdot)$}, if $\Lambda(t_1)\cap\Lambda(0)\neq 0$. 
The dimension of the space $\Lambda(t_1)\cap\Lambda(0)$ is 
called the {\sl multiplicity} of the focal time $t_1$. 
Denote by $\#{\rm foc}_0 
\Lambda(\cdot)\bigl|_{I}\bigr.$ the number of focal times to $0$ on the subset 
$I$ w.r.t. $\Lambda(\cdot)$, counted with their 
multiplicity. If the curve $\Lambda(\cdot)$ is monotone 
increasing, then $\#{\rm foc}_0 
\Lambda(\cdot)\bigl|_{I}\bigr.$ is finite and one can write
\begin{equation}
\label{focn} \#{\rm foc }_0 \Lambda(\cdot)\bigl|_{I}\bigr.= 
\sum_{t\in I} \dim \bigl(\Lambda(t)\cap\Lambda(0)\bigr). 
\end{equation}
 Fix some tuple $\ell=(l_1\ldots,l_s)$ of $s$ linearly independent vectors 
 in 
$\Sigma$, satisfying (\ref{isot}). The time $t_1$ is called 
{\sl focal to the time $0$ w.r.t. the $\ell$-reduction 
$\Lambda(\cdot)^{\ell}$ of the curve $\Lambda(\cdot)$}, if 
$\Lambda^\ell(t_1)\cap\Lambda^\ell(0)\neq {\rm span}\,\ell$ 
or, equivalently, $t_1$ is the focal time to $0$ w.r.t. the 
curve $\overline{\Lambda(\cdot)^\ell}$ in the Lagrange 
Grassmannian $L(({\rm span}\,\ell)^\angle\slash {\rm 
span}\, \ell)$. The multiplicity of the focal time $t_1$ to 
$0$ w.r.t. the $\ell$-reduction $\Lambda(\cdot)^{\ell}$ is 
equal by definition to 
\begin{equation}
\label{inters3} \dim\,\bigl( 
\overline{\Lambda(t_1)^\ell}\cap\overline{\Lambda(0)^\ell}\bigr)=\dim\, 
\bigl(\Lambda(t_1)^\ell\cap\Lambda(0)^\ell\bigr)-s. 
\end{equation}
So, the number of focal times to $0$ on the subset $I$ 
w.r.t. the  $\ell$-reduction $\Lambda(\cdot)^\ell$, 
counted with their multiplicity, is equal to $\#{\rm foc}_0 
\overline{\Lambda(\cdot)^\ell}\bigl|_{I}\bigr.$. 

 It is not 
hard to see that if $\Lambda(\cdot)$ is monotone 
increasing, then $\overline{\Lambda(\cdot)^\ell}$ is 
monotone increasing too. So, the number of focal times to 
$0$ w.r.t. the $\ell$-reduction is also finite. The 
question is what is the relation between the set of points, 
which are focal to $0$ w.r.t. the curves $\Lambda(\cdot)$ 
and its $\ell$-reduction $\Lambda(\cdot)^\ell$? To answer 
this question, we use the fact that if $\Lambda(\cdot)$ is 
monotone curve, then the number $\#{\rm foc } 
\Lambda(\cdot)\bigl|_{(0, T]}\bigr.$ can be represented as 
the intersection index of this curve with a certain 
cooriented hypersurface in $L(\Sigma)$. The advantage of 
this representation is that the intersection index is 
a homotopic invariant. 
 
More precisely, for given Lagrangian subspace $\Lambda_0$ 
denote by $\Mu_{\Lambda_0}$ the following subset of 
$L(\Sigma)$: $$ 
\Mu_{\Lambda_0}=L(\Sigma)\backslash\Lambda_0^{\tr}= 
\{\Lambda\in L(\Sigma):\;\Lambda\cap\Lambda_0\ne 0\}\,. $$ 
Following \cite{Arn0}, the set $\Mu_{\Lambda_0}$ is called 
the {\sl train} of the Lagrangian subspace $\Lambda_0$. 
The set $\Mu_{\Lambda_0}$ is a hypersurface in $L(\Sigma)$ 
with singularities, consisting of the Lagrangian subspaces 
$\Lambda$ such that $\dim\, (\Lambda\cap \Lambda_0)\geq 2$. 
The set of singular points has codimension $3$ in 
$L(\Sigma)$. 
As we have already seen, the tangent space 
$T_{\Lambda}L(\Sigma)$ has a natural identification with 
the space of quadratic forms on $\Lambda$. If $\Lambda$ is 
a non-singular point of the train $\Mu_{\Lambda_0}$, then 
vectors from $T_{\Lambda}L(\Sigma)$ that correspond to 
positive or negative definite quadratic forms are not 
tangent to the train. It defines the canonical 
coorientation of the hyper-surface $\Mu_{\Lambda_0}$ at a 
non-singular point $\Lambda$ by taking as a positive side 
the side of $\Mu_{\Lambda_0}$ containing positive definite 
forms. The defined coorientation permits to define 
correctly the intersection index 
$\Lambda(\cdot)\cdot\Mu_{\Lambda_0}$ 
of an arbitrary continuous curve in the Lagrangian Grassmannian 
$\Lambda(\cdot)$, having endpoints outside 
$\Mu_{\Lambda_0}$:
If $\Lambda(\cdot)$ is smooth and transversally 
intersecting $\Mu_{\Lambda_0}$ in non-singular points, 
then, as usual, every intersection point $\Lambda(\bar t)$ 
with $\Mu_{\Lambda_0}$ adds $+1$ or $-1$ into the value of 
the intersection index according to the direction of the 
vector $\dot\Lambda(\bar t)$ respectively to the positive 
or negative side of $\Mu_{\Lambda_0}$. Further, an 
arbitrary continuous curve $\Lambda(\cdot)$ with endpoints 
outside $\Mu_{\Lambda_0}$ can be (homotopically) perturbed 
to a curve which is smooth and transversally intersects 
$\Mu_{\Lambda_0}$ in non-singular points. Since the set of 
singular points of $\Mu_{\Lambda_0}$ has codimension $3$ in 
$L(\Sigma)$, any two curves, obtained by such perturbation, 
can be deformed one to another by homotopy, which avoids 
the singularities of $\Mu_{\Lambda_0}$. Hence the 
intersection index of the curve, obtained by the 
perturbation, does not depend on the perturbation and can 
be taken as the intersection index of the original curve. 
The intersection index $\Lambda(\cdot)\cdot\Mu_{\Lambda_0}$ 
can be calculated using the notion of the Maslov index of 
the triple of Lagrangian subspaces (see \cite{Agr4}, \cite 
{Arn0}) for details). This implies in particular that if 
the curve $\Lambda:[0,T]\mapsto L(\Sigma)$ is monotone 
increasing with endpoints outside $\Mu_{\Lambda_0}$, then 
\begin{equation}
\label{monotone} 
\Lambda(\cdot)\cdot\Mu_{\Lambda_0}=
\sum_{0\leq t\leq T} \dim \bigl(\Lambda(t)\cap\Lambda_0\bigr)
\end{equation}
Now we are ready to formulate the main result of this 
section: 

\begin{theorem}\label{Mindex}
Let $\Lambda:[0,T]\mapsto L(\Sigma)$ be a monotone 
increasing curve and $\ell=(l_1\ldots,l_s)$ be a tuple of 
$s$ linearly independent vectors in $\Sigma$, satisfying 
(\ref{isot}) and (\ref{inters1}) at $0$. Then on the set 
$(0,T]$ the difference between the number of focal times to 
$0$ w.r.t. the $\ell$-reduction $\Lambda(\cdot)^\ell$ and 
the number of focal times to $0$ w.r.t. the curve 
$\Lambda(\cdot)$ itself , counted with their multiplicity, 
is nonnegative and does not exceed $s$, namely 

\begin{equation}
\label{masineq}
0 \le \#{\rm foc }_0\overline{\Lambda(\cdot)^\ell}\bigl|_{(0,T]}\bigr.-\#{\rm foc }_0
\Lambda(\cdot)\bigl|_{(0,T]}\bigr.\le s\,.
\end{equation}
\end{theorem}
 
{\bf Proof.} 
Since the curves $\Lambda(\cdot)$ and 
$\overline{\Lambda(\cdot)^\ell}$ are monotone increasing, 
for sufficiently small $\varepsilon>0$ the set 
$(0,\varepsilon]$ does not contain the times focal to 0 
w.r.t. both of these curves. Also, without loss of 
generality, one can assume that $T$ is not focal to $0$ 
w.r.t. both of these curves (otherwise, one can extend 
$\Lambda(\cdot)$ as a monotone increasing curve to a 
slightly bigger interval $[0, T+\tilde\varepsilon]$ such 
that $\Lambda(T+\tilde\varepsilon)\cap\Lambda(0)=0$ and 
$\overline{\Lambda(T+\tilde\varepsilon)^\ell}\cap 
\overline{\Lambda(0)^\ell}=0$). So, by relations 
(\ref{focn}), (\ref{monotone}), we have first that 
\begin{equation}
\label{monotone1} \#{\rm foc }_0 
\Lambda(\cdot)\bigl|_{(0,T]}\bigr.= 
\Lambda(\cdot)\bigl|_{[\varepsilon , 
T]}\bigr.\cdot\Mu_{\Lambda(0)}. 
\end{equation}
Besides, from (\ref{inters1}) it follows that 
\begin{equation}
\label{inters4} 
\dim\,\bigl(\Lambda(t)^\ell\cap\Lambda(0)^\ell\bigr)-s= 
\dim\,\bigl(\Lambda^\ell(t)\cap\Lambda(0)\bigr). 
\end{equation} 
Hence, combining (\ref{focn}) and (\ref{monotone}) with 
(\ref{inters3}) and (\ref{inters4}), we get

\begin{equation}
\label{monotone2}
 \#{\rm foc }_0\overline{\Lambda(\cdot)^\ell}\bigl|_{(0,T]}\bigr.=
\Lambda(\cdot)^\ell\bigl|_{[\varepsilon , T]}\bigr.\cdot\Mu_{\Lambda(0)}
\end{equation}

Now we prove the theorem in the case $s=1$. In this case 
$\ell=l_1$. We use the invariance of the defined 
intersection index under homotopies, preserving the 
endpoints. 

Let $a_1(t)$ be as in (\ref{aux}). Denote
\begin{equation}
\label{homot} F(\tau, t)={\rm span}\bigl(\Lambda(t)\cap 
l_1^\angle, (1-\tau)a_1(t)+\tau l_1\bigr) 
\end{equation}  
Note that all subspaces $F(\tau, t)$ are Lagrangian. Let 
$\Phi_\tau:[0,T]\mapsto L(\Sigma)$ and 
$\Gamma_t:[0,1]\mapsto L(\Sigma)$ 
 be the curves, 
satisfying \begin{equation} \label{FG} 
\begin{split} &\Phi_\tau(\cdot)=F(\tau, \cdot),\quad 0\leq\tau\leq 
1;
\\
&\Gamma_t(\cdot)=F(\cdot,t), \quad 0\leq t\leq T. 
\end{split} 
\end{equation}  Then $\Phi_0 (\cdot)=\Lambda(\cdot)$, $\Phi_1( 
\cdot)=\Lambda^{l_1}(\cdot)$, and the curves 
$$\Gamma_\varepsilon(\cdot)\big|_{[0,\tau]}\cup 
\Phi_\tau(\cdot)\big|_{[\varepsilon,T]}\cup 
\Bigl(-\Gamma_T(\cdot)\big|_{[0,\tau]}\Bigr)$$ define the 
homotopy between $\Lambda(\cdot)\big|_{[\varepsilon,T]}$ 
and 
$\Gamma_\varepsilon(\cdot)\cup\Lambda(\cdot)^{l_1}\big|_{[\varepsilon,T]} 
\cup \bigl(-\Gamma_T(\cdot)\bigr)$, preserving the 
endpoints (here $-\gamma(\cdot)$ means the curve, obtained 
from a curve $\gamma(\cdot)$ by inverting the orientation). 
Therefore, $$ 
\Lambda(\cdot)\big|_{[\varepsilon,T]}\cdot\Mu_{\Lambda(0)}= 
\Gamma_\varepsilon(\cdot)\cdot\Mu_{\Lambda(0)}+ 
\Lambda(\cdot)^{l_1}\big|_{[\varepsilon,T]}\cdot\Mu_{\Lambda(0)}- 
\Gamma_T(\cdot)\cdot\Mu_{\Lambda(0)} $$ Using 
(\ref{monotone1}) and (\ref{monotone2}), the last relation 
can be rewritten in the following form 
 \begin{equation}
 \label{masineq1}
  \#{\rm foc }_0 
\overline{\Lambda(\cdot)^{l_1}}\bigl|_{(0,T]}\bigr. -\#{\rm 
foc }_0 \Lambda(\cdot)\bigl|_{(0,T]}\bigr.= 
\Gamma_T(\cdot)\cdot\Mu_{\Lambda(0)}- 
\Gamma_\varepsilon(\cdot)\cdot\Mu_{\Lambda(0)}. 
\end{equation}
So, in order to prove the theorem in the considered case it 
is sufficient to prove the following two relations: 
\begin{eqnarray} 
&~&\label{rT} 0\leq 
\Gamma_T(\cdot)\cdot\Mu_{\Lambda(0)}\leq 1, \\ &~& 
\label{reps} \exists \varepsilon_0>0\,\,{\rm s.}\,{\rm 
t.}\,\, \forall \varepsilon_0\geq \varepsilon>0:\quad 
\Gamma_\varepsilon(\cdot)\cdot\Mu_{\Lambda(0)}=0. 
\end{eqnarray}
%for sufficiently small $\varepsilon>0$.

{\bf a)} Let us prove (\ref{rT}). If $l_1\in \Lambda(T)$, 
then by definition $\Gamma_T(\tau)\equiv \Lambda(T)$. 
Since, by our assumptions, $\Lambda(T)\cap \Lambda(0)=0$, 
we obviously have $\Gamma_T(\cdot)\cdot\Mu_{\Lambda(0)}=0$. 

If $l_1\not\in \Lambda(T)$, then $\dim\, 
(\Lambda(0)+\Lambda(T)\cap(l_1)^\angle)=2n-1$. In 
particular, it implies that 
\begin{equation}
\label{dim1}0\leq  
\dim\,\bigl(\Gamma_T(\tau)\cap\Lambda(0)\bigr)\leq 1. 
\end{equation} 
Further, let $p:\Sigma\mapsto 
\Sigma/(\Lambda(0)+\Lambda(T)\cap(l_1)^\angle)$ be the 
canonical projection on the factor space. Then from 
(\ref{homot}) and (\ref{FG}), using standard arguments of 
Linear Algebra, it follows that 
$\Gamma_T(\tau)\cap\Lambda(0)\neq 0$ if and only if 
\begin{equation}
\label{pconv} (1-\tau) p\bigl(a_1(T)\bigr)+\tau p(l_1)=0 
\end{equation}
Since, by assumptions, $\Gamma_T(0)\cap \Lambda(0)=0$ 
(recall that $\Gamma_T(0)=\Lambda(T)$), the equation 
(\ref{pconv}) has at most one solution on the segment 
$[0,1]$. In other words, the curve $\Gamma_T(\cdot)$ 
intersects the train $\Mu_{\Lambda(0)}$ at most ones and 
according to (\ref{dim1}) the point of intersection is 
non-singular.

 Finally, the curve $\Gamma_T(\cdot)$ is 
monotone non-decreasing, i.e. its velocities 
$\frac{d}{d\tau}\Gamma_T(\tau)$ are non-negative definite 
quadratic forms for any $\tau$. Indeed, since 
$\Lambda(T)\cap (l_1)^{\angle}$ is the common space for all 
$\Gamma_T(\tau)$, one has 
$\frac{d}{d\tau}\Gamma_T(\tau)\big|_{\Lambda(T)\cap 
(l_1)^{\angle}} \equiv 0$. On the other hand, if we denote 
by $c(\tau)=(1-\tau)a_1(t)+\tau l_1$, then by 
(\ref{tangdef}), one has 
$$\frac{d}{d\tau}\Gamma_T(\tau)\bigl(c(\tau)\bigr)= 
\sigma\bigl(c'(\tau),c(\tau)\bigr)=\sigma\bigl(l_1-a_1(T),(1-\tau) 
a_1(T)+\tau l_1\bigr)=\sigma(l_1, a_1)>0$$ (the last 
inequality follows from (\ref{ups}), (\ref{sigal}) and the 
assumption about monotonicity of $\Lambda(\cdot)$). 

So, $\frac{d}{d\tau}\Gamma_T(\tau)$ are non-negative 
definite quadratic forms. Hence in the uniquely possible 
point of intersection of $\Gamma_T(\cdot)$ with the train 
$\Mu_{\Lambda(0)}$ the intersection index becomes equal to 
$1$. This proves (\ref{rT}). 

{\bf b)} Let us prove (\ref{reps}). 
 Take a Lagrangian subspace 
$\Delta$ such that $l_1\in\Delta$ and 
$\Delta\cap\Lambda(0)=0$. Then there exists $\varepsilon_0$ 
such that 

\begin{equation}
\label{pitchfork1} 
\Lambda(\cdot)\big|_{[0,\varepsilon_0]}\subset \Delta 
^\pitchfork. 
\end{equation}
Similarly to the arguments in {\bf a)}, for any 
$0<\varepsilon\leq \varepsilon_0$ the curve 
$\Gamma_{\varepsilon}(\cdot)$ intersects the train 
$\Mu_\Delta$ once. But by construction this unique 
intersection occurs at $\tau=1$. Indeed, 
$\Gamma_\varepsilon(1)=\Lambda(\varepsilon)^{l_1}$, hence 
$l_1\in\Gamma _\varepsilon(1)\cap\Delta$. In other words, 
\begin{equation}
\label{pitchfork2} \Gamma 
_{\varepsilon}(\cdot)\big|_{[0,1)} 
\subset\Delta^\pitchfork. 
\end{equation}
Further, one can choose a symplectic basis in $\Sigma$ such 
that $\Sigma=\mathbb R^n\times \mathbb R^n$, the symplectic 
form $\sigma$ is as in (\ref{sympf}),  
$\Lambda(0)=0\times\mathbb R^n$, 
and $\Delta=\mathbb R^n\times 0$. 
By (\ref{pitchfork1}) and (\ref{pitchfork2}), there exists 
two one parametric families of symmetric matrices $S_t$, 
$0\leq t\leq \varepsilon_0$ and $C_\tau$ ,$0\leq \tau<1$ 
such that $\Lambda(t)=\{(S_t p,p):p\in {\mathbb R}^n\}$ and 
$\Gamma_\varepsilon (\tau)=\{(C_\tau p,p):p\in {\mathbb 
R}^n\}$. Since the curve $\Lambda(\cdot)$ is monotone 
increasing and the curve $\Gamma_\varepsilon(\cdot)$ is 
monotone nondecreasing, for any $0\leq \tau <1$ the 
quadratic forms $p\mapsto \langle C_\tau p,p\rangle$ are 
positive definite, while $S_0=0$. It implies that 
\begin{equation}
\label{tauf} \forall 0\leq \tau 
<1:\quad \Gamma_\epsilon (\tau)\cap\Lambda(0)=0.
\end{equation}
Note also that for sufficiently small $\varepsilon>0$ 
\begin{equation} \label{1f} 
\Gamma_\varepsilon(1)\cap\Lambda(0)=0. 
\end{equation}
 Indeed, 
$\Gamma_\varepsilon(1)=\Lambda^{l_1}(\varepsilon)$  and a  
sufficiently small $\varepsilon>0$ is not a focal time for 
the $l_1$-reduction $\Lambda^{l_1}(\cdot)$, which 
according to (\ref{inters4}) is equivalent to the fact that 
$\Lambda^{l_1}(\varepsilon)\cap\Lambda(0)=0$ and hence to 
(\ref{1f}). By (\ref{tauf}) and (\ref{1f}), for 
sufficiently small $\varepsilon>0$ the curve 
$\Gamma_\epsilon(\cdot)$ does not intersect the train 
$\Mu_{\Lambda(0)}$. The relation (\ref{reps}) is proved, 
which completes the proof of our theorem in the case $s=1$. 

The case of arbitrary $s$ can be obtained immediately from 
the case $s=1$ by induction, using the fact that 
\begin{equation} \label{redind} 
\Lambda(\cdot)^{(l_1,\ldots, l_s)}= 
\left(\Lambda(\cdot)^{(l_1,\ldots , l_{s-1})}\right)^{l_s}.\quad  \Box 
\end{equation} 

\begin{remark} {\rm Note that in the case $s=1$ from Theorem 
\ref{Mindex} it follows immediately that the sets of focal 
times ( to $0$) w.r.t monotone increasing curve and its 
reduction are alternating. 
Also, for any $s$ the first focal time to $0$  w.r.t. the 
reduction does  exceed the first focal time w.r.t. the curve 
itself.} $\Box$ 
\end{remark}

All constructions above are directly related to the notion 
of focal points of a dynamical Lagrangian distributions and 
$(\vec{\mathcal H},{\mathcal D})$ and its reduction by a 
tuple $\mathcal{G}=(g_1,\ldots, g_s)$ of $s$ involutive 
first integrals, defined in Introduction. Note that the 
point $\lambda_1=e^{t_1\vec{\mathcal H}}\lambda_0$ is focal 
to $\lambda_0$ w.r.t. the pair $(\vec{\mathcal H},{\mathcal 
D})$ along the integral curve $t\mapsto e^{t\vec{\mathcal 
H}}\lambda_0$ of $\vec{\mathcal H}$ if and only if the time 
$t_1$ is focal to $0$ w.r.t the Jacobi curve 
$J_{\lambda_0}(\cdot)$ attached at the point $\lambda_0$, 
while $\lambda_1=e^{t_1\vec{\mathcal H}}\lambda_0$ is focal 
to $\lambda_0$ w.r.t. the $\mathcal {G}$-reduction of the 
pair $(\vec{\mathcal H},{\mathcal D})$ along the integral 
curve $t\mapsto e^{t\vec{\mathcal H}}\lambda_0$ of 
$\vec{\mathcal H}$ if and only if $t_1$ is focal to $0$ 
w.r.t $\bigl(\vec g_1(\lambda_0),\ldots,\vec 
g_s(\lambda_0)\bigr) $-reductions of the Jacobi curves 
$J_{\lambda_0}(\cdot)$ attached at $\lambda_0$. Translating 
Theorem \ref{Mindex} into the terms of dynamical Lagrangian 
distribution, we have immediately the following 
\begin{corollary}
\label{Mcor} 
 If the dynamical Lagrangian 
distribution $(\vec{\mathcal H},{\mathcal D})$ is monotone 
increasing and the Hamiltonian ${\mathcal H}$ admits a 
tuple $\mathcal{G}=(g_1,\ldots, g_s)$ of $s$ involutive 
first integrals satisfying (\ref{inters2}), then along any 
segment of the integral curve of $\vec{\mathcal H}$ the 
difference between the number of the focal points to the 
starting point of the segment w.r.t. the ${\mathcal G} 
$-reduction of the pair $(\vec{\mathcal H},{\mathcal D})$ 
and the number of the focal points to the starting point of 
the segment w.r.t. the pair $(\vec{\mathcal H},{\mathcal 
D})$ itself , counted with their multiplicity \footnote{ 
Here we do not count the starting point of the segment as 
the focal point to itself.}, is nonnegative and does not 
exceed $s$. 
\end{corollary}

\begin{ex}
{\sl (Plane 3-body problem with equal masses: focal points 
of the 8-shaped orbit)} {\rm 
%To conclude our analysis of the relation between the 
%reduction and the behavior of the Jacobi curve we consider 
%one more example that illustrates the Theorem \ref{Mindex}. 
The following example illustrates the Theorem \ref{Mindex}. 
In 2000, A. Chenciner and R.Montgomery proved the existence 
of a new periodic solution of the planar 3-body problem 
with equal masses - the $8$-shaped orbit or just {\sl the 
Eight} ~\cite{C&M}. In the plane of the motion each body 
moves along the same $8$-shaped orbit, symmetric w.r.t. the 
point of its self-intersection, coinciding with the center 
of mass of the bodies. The configuration space is $M=\real 
^6$. As an initial point in the phase space $T^*M$ we take 
the point $\lambda_0$ such that its projection on the 
configuration space is a collinear configuration, i.e. one 
of the bodies lies in the middle of the segment, connecting 
the other two. 
%This orbit was found analytically by minimizing the action 
%functional over the space of loops that satisfy some 
%properly chosen symmetry conditions, so that the 
%$\frac{1}{12}$-th piece of the Eight provides the minimum 
%of the action. 

In ~\cite{Chtch} there were found the  
focal points to $\lambda_0$ along  the Eight w.r.t. the $g$-reduction of the 
Lagrangian dynamical distribution $\bigl(\Ham,\Pi(M)^g\bigr)$,
where $H$, $g$ are  as in Example 9. In particular it was shown 
numerically that the 8-shaped orbit contains three such focal 
points along its period $T$, and the first focal time 
$\tau_1\approx0.52T$. 
      
Let  $e^{t_i\Ham}\lambda_0$ be the $i$th focal point w.r.t. $\bigl(\Ham,\Pi(M)\bigr)$ along the Eight, and 
let  $e^{\tau_i\Ham}\lambda_0$ be the $i$th focal point w.r.t. its $g$-reduction $\bigl(\Ham,\Pi(M)^g\bigr)$ along the same curve.
In the following table we present the result 
of the numerical computation of $t_i$ and $\tau_i$ on the interval $(0,3T]$
(all this focal points have the multiplicity 1): 

$$\begin{array}{|c|c|c|c|c|c|c|c|c|c|c|c|} \hline 
i&1&2&3&4&5&6&7&8&9&10&11\\ \hline 
\tau_i/T\approx&0.52&0.76&0.95&1.08&1.52&1.56&1.88&2.05&2.29&2.49&2.65\\ 
\hline 
t_i/T\approx&0.76&0.95&1.08&1.42&1.54&1.88&2.05&2.28&2.45&2.65&\,\\ 
\hline 
\end{array}
$$
We observe  that $\tau_i\leq t_i\leq \tau_{i+1}$, $1\leq i\leq 10$, as  was expected by Theorem \ref{Mindex}. 
$\Box$}
\end{ex}

\end{document}